\documentclass[10pt,a4paper]{article}
\usepackage{makeidx}
\usepackage{amssymb}
\usepackage{amsmath}
\usepackage{graphicx}
\usepackage{graphics}
\usepackage{times}
\usepackage[T1]{fontenc}
\usepackage[thinspace,thinqspace,squaren]{SIunits}
\usepackage{graphicx}
\usepackage{amsthm}
\usepackage{color}
\usepackage{hyperref}
\usepackage{multirow}

\usepackage{ctable}
\allowdisplaybreaks

\usepackage{algorithm}
\usepackage{algorithmic}

\twocolumn

\begin{document}
\date{}

\title{On the Take-off of Airborne Wind Energy Systems\\ Based on Rigid Wings}


\author{L. Fagiano and S. Schnez%
\thanks{This is the pre-print of a paper submitted for possible publication to the Elsevier journal  \emph{Energy}. The authors are with ABB Switzerland Ltd., Corporate Research, 5405 Baden-D\"{a}ttwil - Switzerland. E-mail addresses: \{ lorenzo.fagiano $|$ stephan.schnez\}@ch.abb.com. Both authors contributed equally to this publication.}}
\maketitle

\begin{abstract}
The problem of launching a tethered aircraft to be used for airborne wind energy  generation is investigated. Exploiting well-assessed physical principles, an analysis of three different take-off approaches is carried out. The approaches are then compared on the basis of quantitative and qualitative criteria introduced  to assess their technical and economic viability. Finally, a deeper study of the concept that is deemed the most viable one, i.e. a linear take-off maneuver combined with on-board propellers, is performed by means of numerical simulations. The latter are used to refine the initial analysis in terms of power required for take-off, and further confirm the viability of the approach.
\end{abstract}



\section{Introduction}
The term airborne wind energy (AWE) refers to a class of wind power
generators that exploit tethered aircrafts to convert wind energy
 into electricity \cite{AWEbook,FaMi12}. The benefits
of AWE systems, compared to traditional wind turbines, are essentially two: lower construction and
installation costs and the possibility to reach higher altitudes,
where faster and steadier winds blow. According to the current
estimates the combination of these two benefits should render AWE systems
competitive with the established energy sources, including fossil
fuels \cite{FaMP09}, in terms of both cost of energy and land occupation. The first papers and patents concerned with
AWE appeared in the late 1970s (see e.g. \cite{Mana76,Loyd80}),
yet only in recent years a significant and growing research
effort has been undertaken by both small companies and universities
to develop such concepts via theoretical, numerical and experimental
methods \cite{AWEbook}. AWE is still in its infancy and no
commercial 
system exists; 
however, thanks to the continuous progresses
that are being achieved, a relatively well-established set of
few different approaches has emerged, while other, less promising
ideas have been abandoned.

Today, AWE systems can
be classified by the way the lift force that keeps the aircraft
airborne is generated -- either aerodynamic lift \cite{ch20-ZiHa15,ch21-MiTM14,ch23-vdVlPS14,ch24-BRKGS14,ch26-RuSo14,ch28-Vand14}, or aerostatic lift \cite{ch30-VeGR14} -- and by the placement of the
electrical generators - either on-board of the aircraft
\cite{ch28-Vand14,ch30-VeGR14} or on the ground
\cite{ch20-ZiHa15,ch21-MiTM14,ch24-BRKGS14,ch23-vdVlPS14,ch26-RuSo14}.
Among the systems that exploit aerodynamic lift and ground-level
generators, a further distinction can be made between concepts that
rely on rigid wings \cite{ch26-RuSo14}, similar to gliders, and
concepts that employ flexible wings, like power kites
\cite{ch20-ZiHa15,ch21-MiTM14,ch23-vdVlPS14,ch24-BRKGS14}.
Small-scale prototypes (10-50 kW of rated power) of all of the
mentioned concepts have been realized and successfully tested to
demonstrate their power generation functionalities. Moreover,
scientific contributions concerned with several different technical
aspects, primarily aerodynamics
\cite{BSTR14,ch16-BrSO14,ch17-BSTR14,ch18-GoLu14,ch19-LRBLJP14}
and controls
\cite{HoDi07,CaFM09c,BaOc12,ErSt12,FZMK14,FHBK14,FeSc14,ch12-ZaGD14,ErSt14,ZaFM15}
but also resource assessment \cite{ArCa09,ch5-Arch15}, economics
\cite{FaMP09,ch7-ZiHa15}, prototype design \cite{FaMa15} and
power conversion \cite{SHVDD15}, have recently appeared, gradually
improving and expanding our understanding of such systems.

Despite the steady and promising development of the field, the complexity
and multidisciplinary nature of AWE systems are such that several
relevant aspects still need to be addressed in order to ultimately
prove the technical and economic feasibility of the idea. One of such aspects is
the take-off of the aircraft, particularly for concepts that employ
rigid wings and ground-level generation. In fact, while systems
with on-board generation \cite{ch28-Vand14, ch30-VeGR14}, as well as kite-based
systems with ground generation \cite{ch20-ZiHa15} are able to take-off autonomously from a quite compact ground
area, the same functionality for AWE systems with rigid wings and
ground-level generators has not been achieved yet. There is
evidence of autonomous take-off of this class of generators
\cite{ampyx}; however by using a winch launch that requires
a significant space in all directions in
order to adapt to the wind conditions during take-off. As a
consequence, one of the main advantages of AWE systems, i.e. the
possibility of being installed in a large variety of
locations at low costs, 
might be lost due to the need of a large area of land
suitable for the take-off. So far, this issue has been addressed only to a limited extent within the scientific community. 
In Ref.~\cite{ZaGD13}, a rotational take-off is studied and simulated; however the focus is on the control and optimization aspects of this approach, rather than on its economic viability and the comparison with other possible methods. In Ref.~\cite{Bont10}, an analysis of several approaches is first carried out, considering different performance criteria, and three alternatives are deemed the most promising: buoyant systems, linear ground acceleration plus on-board propeller, and rotational take-off. Then, the rotational take-off is examined in more detail by means of numerical simulations.

In order to contribute to address this important problem, we present here an analysis of three
candidate approaches to realize the take-off of a rigid tethered
aircraft with ground-based generation. This is the concept which is also pursued by the company Ampyx Power \cite{ampyx,ch26-RuSo14}. More specifically, we
compare a vertical lift approach with on-board vertical-axis propellers, a rotational
take-off like the one considered in Refs.~\cite{ZaGD13,Bont10}, and a linear take-off
technique combined with on-board horizontal-axis propellers. The analysis is
instrumental to carry out a comparison among the considered
approaches, based on a series of performance criteria that we
introduce in order to quantify their viability. The analysis and the
subsequent comparison represent the first main contribution that the
present paper adds to the existing scientific literature. Then, we study in more depth the concept that
is deemed the most viable one, i.e. the
linear take-off maneuver combined with on-board
propellers. In particular, we derive a dynamical model of the system
that includes realistic aerodynamic coefficients, as well as
friction and inertia, and we use it to refine the initial
analysis in terms of power required for take-off.
Since the system is unstable in open-loop, we also develop the
feedback control algorithms required to stabilize the take-off
maneuver and carry out the numerical simulations.

The paper is organized as follows: section \ref{S:preliminaries}
provides more details on the considered type of AWE system, which
are needed to formulate the problem that we address in a rigorous way, and a brief description of the considered take-off approaches. The performance criteria are given in section \ref{S:preliminaries}, too. section \ref{S:evaluation} presents the analysis of the three take-off concepts using basic physical equations. The numerical simulation study is reported in section
\ref{S:simulation}. Final conclusions are drawn in section
\ref{S:conclusions} together with a discussion of future
developments of this research.

\section{Preliminaries and problem formulation}\label{S:preliminaries}
We first describe the system under consideration and introduce the physical equations that link the
main lumped design parameters to the generated mechanical power. These equations can be employed in a first-approximation dimensioning phase of the AWE generator and are used here to compute one of our performance criteria. For the complete details and derivation of the
equations recalled in the following, we refer to \cite{Loyd80,FaMP11,FaMi12, AWEbook}.

\subsection{Airborne wind energy systems based on rigid aircrafts and ground-level generation}\label{SS:system_description}

\begin{figure}[!htb]
 \begin{center}
  \includegraphics[trim= 0cm 0cm 0cm 0cm,width=8cm]{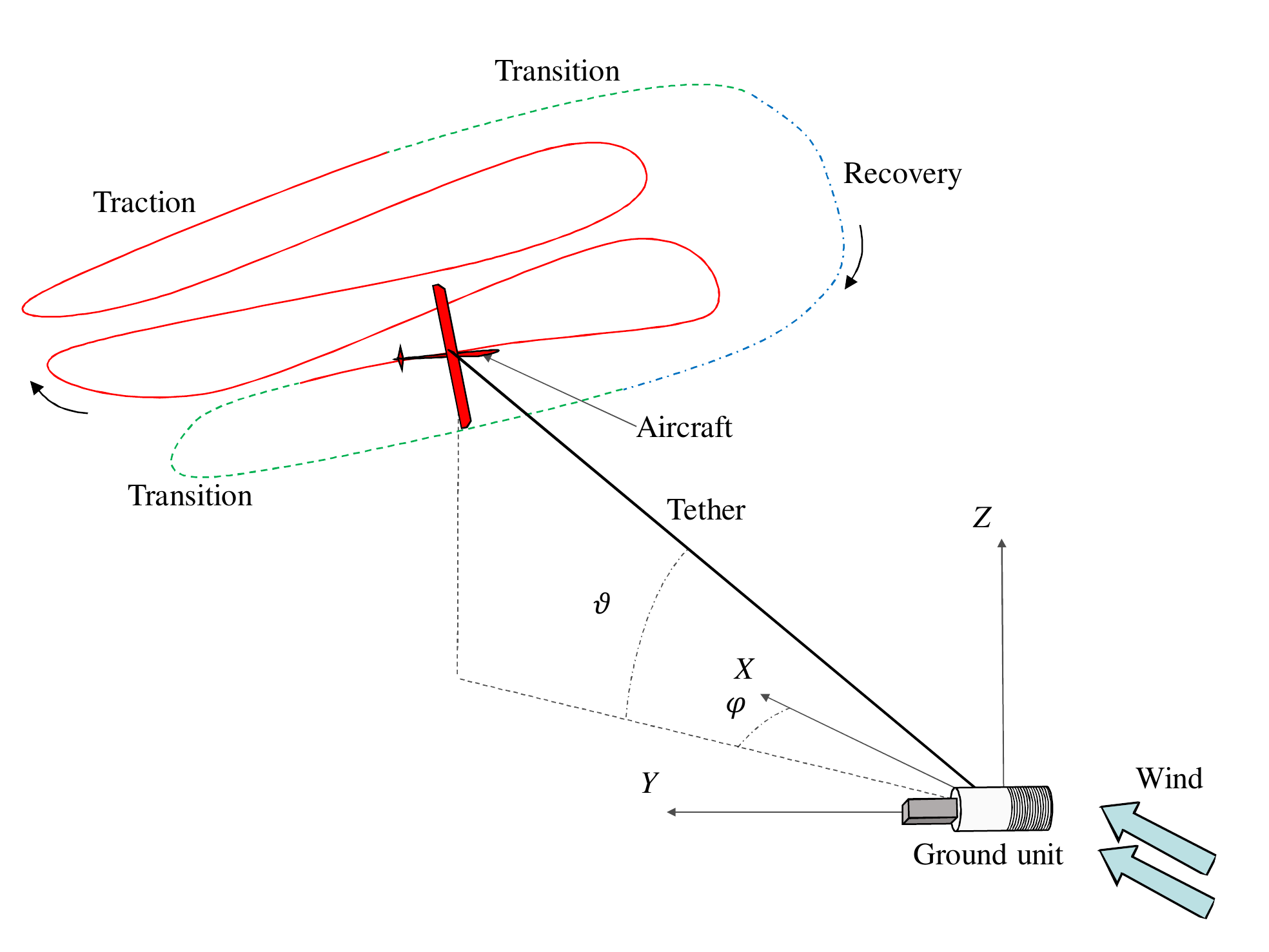}
  \caption{Sketch of the considered AWE generator and its working principle during power
  production. In the traction phase (red solid line) the aircraft is controlled to
  follow figure-of-eight patterns in crosswind conditions, and the
  tether is reeled out under large load from the drum installed in the GU.
  In the retraction phase (blue dash-dotted line), the aircraft is controlled to glide towards
  the ground station, and the tether is reeled-in under small load.
  Two transitions (green dashed lines) link the traction and
  retraction phases. The aircraft position with respect to the incoming wind can be defined by the
  elevation angle $\vartheta$ and the azimuth angle $\varphi$.}
  \label{F:system_sketch}
 \end{center}
\end{figure}

The considered AWE system is composed of a rigid
aircraft, a ground unit (GU), and a tether connecting them, as depicted
in Figure \ref{F:system_sketch}. The aircraft is equipped with
sensors, actuators and on-board intelligence to attain autonomous
flight and realize the flight patterns required to generate power,
as well as with communication capabilities to exchange information
with the GU and possibly with other systems and infrastructure
nearby. 

The GU consists of several subsystems, the main
ones being a drum, around which the tether is coiled, an electric
machine (generator/motor), linked to the drum through a mechanical
transmission system, and the power electronic system to control the
generator and to convert mechanical power into electrical one and
vice-versa. 

The described AWE system generates energy by means of a cyclic operating principle composed essentially by four phases: the power generation (or traction) phase, the retraction phase, and two transition phases linking them, shown in Figure \ref{F:system_sketch}. During the traction phase, the on-board control system steers the aircraft into figure-of-eight patterns under crosswind conditions. The generated aerodynamic forces exert a large traction load on the line, which is reeled-out from the drum. The electric machine exerts a torque on the drum in order to achieve a desired reel-out speed and to produce power. In particular, an aircraft with effective area $A$, aerodynamic lift and drag coefficients $C_l$ and $C_d$, respectively, flying at a relative elevation $\vartheta$ and azimuth $\varphi$ with respect to a wind flow of speed $W$ (see Figure \ref{F:system_sketch}), exerts a traction load $T$ on the tether approximately equal to \cite{Loyd80,FaMP11,FaMi12}:
\begin{equation}\label{E:tether_force}
\begin{array}{rcl}
T(t)&\simeq&\frac{1}{2}\rho A \dfrac{C_l(t)^3}{C_{d,eq}(t)^2}\\&& \left(W(t)\cos{(\varphi(t))}\cos{(\vartheta(t))}-\dot{l}(t)\right)^2
\end{array}
\end{equation}
where $t$ is the continuous time variable, $\rho$ is the air density, $C_{d,eq}\doteq C_d(t)+\frac{d_l\,l(t)\,C_{d,l}}{4\,A}$ is the equivalent drag coefficient (taking into account the drag of both the aircraft and the line), $l$ is the length of the line, assumed straight, $d_l$ its diameter, $C_{d,l}$ its drag coefficient, and $\dot{l}\doteq\frac{dl}{dt}$ is the tether reeling speed. For $\dot{l}>0$, the line is reeled out from the drum, hence effectively decreasing the apparent wind speed parallel to the tether direction, given by $W\cos{(\varphi)}\cos{(\vartheta)}$. The tether force $T(t)$ multiplied with the reeling speed $\dot{l}(t)$ provides an estimate of the instantaneous mechanical power $P_m(t)$ generated during the traction phase:
\begin{equation}\label{E:power}
	P_m(t)\simeq T(t)\,\dot{l}(t).
\end{equation}
The maximum generated power is achieved when the reeling speed is equal to 1/3 of the absolute wind speed projected along the line direction, i.e. $\dot{l}=\dfrac{1}{3}W\cos{(\varphi)}\cos{(\vartheta)}$, and ideally with $\varphi=\vartheta =0$. In this case, the obtained mechanical power is:
\begin{equation}\label{E:power2}
\begin{array}{rcl}
	P_m^*(t)&\simeq& \dfrac{2}{27}\rho\,A\,.\dfrac{C_l(t)^3}{C_{d,eq}(t)^2} W(t)^3.
\end{array}
\end{equation}

For the sake of estimating the generated power, the mass of the airborne components is irrelevant as a first approximation, since the weight and apparent forces of the aircraft and of the line are significantly smaller than the force acting on the tether during the traction phase. On the other hand, this parameter clearly plays a crucial role when discussing take-off approaches. In order to evaluate a given take-off technique on a quantitative basis, the total mass of the aircraft $m$ has to be linked to the system's capability in terms of force and power. Such a link is given by the so-called wing loading $w_l$, i.e. the ratio between $m$ and the effective aerodynamic area $A$:
\begin{equation}\label{E:wing_loading}
m = w_l\,A.
\end{equation}
The total mass of the aircraft is the sum of $m$ and of the additional mass $\Delta m_i$ required for the take-off capability. This will be discussed further in section \ref{performance_criteria}.

%

\subsection{Take-off approaches}\label{SS:approaches}

Here, we briefly describe the three take-off concepts under consideration.

\textbf{Vertical take-off with rotors}. In this approach, the aircraft is equipped with vertical-axis propellers which provide enough lift to take-off vertically. In the framework of ground-level generation, this approach is pursued by the company TwingTec \cite{twingtec}. In the AWE field, the company Makani Power owned by Google \cite{Makani,ch28-Vand14} employs this approach for the take-off and landing their system with on-board power generation.

\textbf{Rotational take-off}. This is the only proposal for rigid-wing systems which has been studied in the literature with numerical simulations in addition to static equations \cite{Bont10,ZaGD13}. In this approach, the hull of the aircraft is initially attached at the tip of a rotating arm. When the tangential speed of the arm is large enough, the aircraft takes off exploiting its aerodynamic lift and the tether is gradually extended out of the rotating arm until a certain altitude is reached. Then, the rotating arm is gradually stopped while the aircraft transitions into power-generating mode. The company EnerKite \cite{enerkite,ch24-BRKGS14} is implementing this concept for its AWE system.

\textbf{Linear take-off with on-board propellers}. In this approach the aircraft is accelerated on a rectilinear path up to take-off speed by an external source of power, for example the winch itself or a linear motion system. 
Horizontal-axis on-board propellers are then employed to sustain the forward speed during the climb to the operational altitude. This approach was briefly analyzed and deemed promising in Ref.~\cite{Bont10}, but without carrying out a deeper analysis by  means of e.g. numerical simulations. The company Ampyx Power \cite{ampyx, ch26-RuSo14} pursues a similar take-off concept as the one discussed here.

In the remainder of this paper, we will use the subscripts $_1,\,_2,\,_3$ respectively for the vertical, rotational and linear take-off approaches described above.

\subsection{Performance criteria and problem formulation}\label{performance_criteria}

A well-established metric to compare different electric power generation schemes on economic grounds is the levelized cost of electricity (LCOE). In our case, additional components or land occupation required to implement the take-off approach will increase upfront costs (and potentially maintenance costs) and will lead to an increase in the LCOE of the AWE system, as compared to the same system without the take-off functionality. Hence, when comparing different take-off approaches, their impact on the LCOE should be assessed. However, the precise calculation of the LCOE is not feasible 
for new power generation concepts like AWE systems.

Rather than the LCOE, we will therefore consider a series of other quantitative and qualitative criteria which are easier to evaluate based on the existing know-how of AWE generators, and which are related to the system's cost, complexity and required land occupation. If a specific take-off approach performs well according to these criteria, we can expect that the impact on the LCOE of the AWE system will be small.

The quantitative criteria are:
\begin{enumerate}
\item[\textbf{C1}] The additional power installed on-board and on the ground, relative to the peak mechanical power of the system, required to carry out the take-off procedure:
    \begin{equation}\label{E:crit_power}
    \begin{array}{l}
    \overline{P}_{g,i}\simeq\eta_{P_g,i}\;P_m^*\\
    \overline{P}_{ob,i}\simeq\eta_{P_{ob},i}\;P_m^*
    \end{array}
    \end{equation}
    where $\overline{P}_g$ and $\overline{P}_{ob}$ stand for the peak ground and on-board power, respectively, and $i=1,\,2,\,3$ refers to the three considered take-off approaches. The higher the values of $\eta_{P_g,i},\,\eta_{P_{ob},i}$, the worse the approach.\\

\item[\textbf{C2}] The additional on-board mass, relative to the aircraft's mass without the system required for the take-off:
    \begin{equation}\label{E:crit_mass}
    \Delta m_i\simeq\eta_{m,i}\;m.
    \end{equation}
    Although, as recalled in section \ref{SS:system_description}, the mass does not impact the maximum power generation in a first approximation, it is an important parameter for the controllability and maneuverability of the system and for its capability to operate in a wide range of wind conditions \cite{FZMK14}. Again, the higher $\eta_{m,i}$, the worse the approach.\\

\item[\textbf{C3}] The ground area occupied by the take-off system, indicated with $A_{g,i}$:
    \begin{equation}\label{E:ground_area}
    A_{g,i}\simeq \underline{A}_{g,i}+\eta_{A_g,i}\;A,
    \end{equation}
    The higher $\underline{A}_{g,i},\,\eta_{A_g,i}$, the worse the approach.
\end{enumerate}

The qualitative criteria that we consider are:
\begin{enumerate}
	\item[\textbf{C4}] The complexity and cost of the apparatus that needs to be added to the system for the take-off functionality.
	\item[\textbf{C5}] The capability to take-off under most wind conditions (including no wind).
\end{enumerate}

The problem we will address in the next section is to carry out a comparison of the three considered approaches in light of criteria C1-C5. In particular, we will derive equations that allow to compute the quantitative criteria C1-C3, and we will assess the criteria C4-C5 on the basis of the knowledge on AWE systems available in the literature and of our own hands-on experience.


\section{Assessment of take-off concepts for rigid-wing AWE systems}\label{S:evaluation}

In the following three sections \ref{SS:rotorcraft}, \ref{SS:rotational}, and \ref{SS:linear}, we will introduce the relevant assumptions and derive the governing equations of the considered take-off approaches. Quantitative results and the related discussion will be presented in section \ref{SS:conclusion}.


\subsection{Vertical take-off with rotors}\label{SS:rotorcraft}

According to the Actuator Disk Theory \cite{Horl78}, the thrust through a disk with area $A_\mathrm{prop}$ is
\begin{equation}\label{E:vertical_thrust}
	T=\frac{1}{2} \rho A_\mathrm{prop} \left(v^2_\mathrm{out}-v^2_\mathrm{in}\right),
\end{equation}
where the velocities are taken far in front and far behind the disk. The associated power is then
\begin{equation}\label{E:vertical_power}
	P_{ob,1}=\frac{1}{2}(v_\mathrm{out}+v_\mathrm{in})\,T.
\end{equation}
In order to lift an object with vertical velocity $v_\mathrm{c}$ and mass $m$, the thrust must equal the weight, $T=m g$. By setting $v_\mathrm{in}=v_\mathrm{c}$ with $v_\mathrm{c}$ being the desired climb velocity, considering a conversion efficiency $\eta<1$ between mechanical power at the shaft and fluid-dynamic power, and solving Eqs.~\eqref{E:vertical_thrust} and \eqref{E:vertical_power} for $P_{ob,1}$, it then follows that the required take-off power is
\begin{equation}\label{heli_power}
	P_{ob,1} = \dfrac{(m+\Delta m_1)g}{\eta}\left(\sqrt{\frac{(m+\Delta m_1) g}{2\rho A_\mathrm{prop}} + \frac{v^2_\mathrm{c}}{4}} + \dfrac{1}{2} v_\mathrm{c}\right).
\end{equation}

In our assessment, for the sake of computing $P_{ob,1}$, we will consider a wing with wingspan $d$ and aspect ratio (i.e. wingspan divided by the chord) $\lambda$, and we will assume that the aircraft employs two propellers with a diameter equal to the chord length, i.e. $d/\lambda$. Thus, we have $A=d^2/\lambda$ and 
$A_\mathrm{prop}=\dfrac{\pi\,d^2}{2\,\lambda^2}$. With regard to the additional on-board mass $\Delta m_1$, this is given mainly by the onboard batteries and electric motors that drive the propellers. The required battery mass is calculated from the energy density of lithium-polymer batteries $E_\text{batt}$ and the required power $P_{ob,1}$, target altitude $h$ and climb speed $v_\mathrm{c}$ (i.e. the climb duration is $h/v_\mathrm{c}$). The power density of an electric motor is indicated by $E_\text{mot}$. The resulting equation for the additional on-board mass is:
\begin{equation}\label{E:heli_deltam}
\Delta m_1= P_{ob,1}\left(\dfrac{h}{v_\mathrm{c}\,E_\text{batt}}+\dfrac{1}{E_\text{mot}}\right)
\end{equation}
We solve the system of Eqs.~\eqref{heli_power} and \eqref{E:heli_deltam} to compute the required take-off power, in order to account also for the additional mass.

Finally, as regards the occupied ground area, we assume that the vertical take-off can be carried out with all possible angles between the wing and the nominal wind speed. Hence, we have
\begin{equation}\label{E:heli_area}
A_{g,1}=\dfrac{\pi d^2}{4}=\dfrac{\pi \lambda}{4}A
\end{equation}


\subsection{Rotational take-off}\label{SS:rotational}

\begin{figure}[t]
	\begin{center}
		\includegraphics[width= \linewidth]{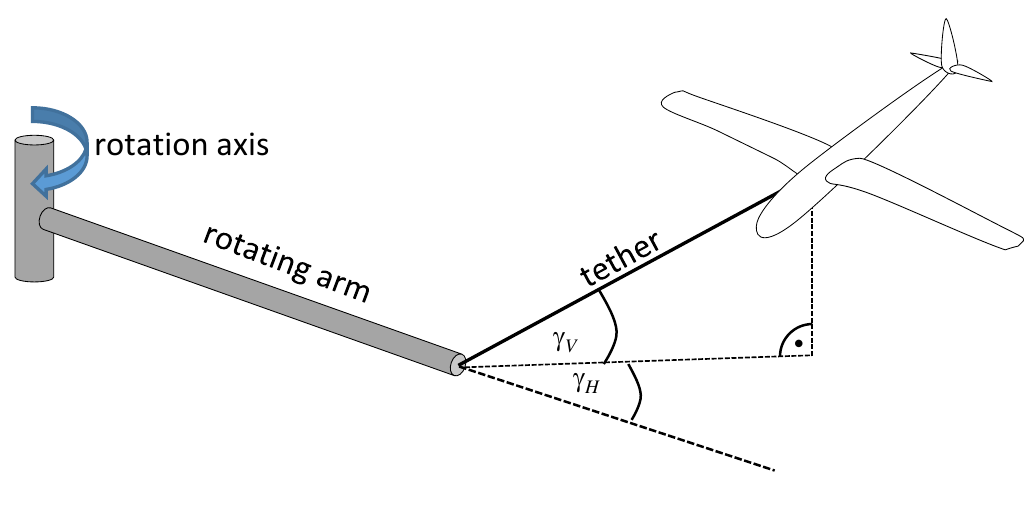}
		\caption{Sketch of an aircraft attached to a rotating arm via the tether during a rotational start. The azimuth of the plane is given by the angle $\gamma_h$; the angle $\gamma_v$ denotes the angle between the tether and the plane of the rotating arm.\label{Fig_RotStart}}
	\end{center}
\end{figure}

\begin{figure}[t]
	\begin{center}
		\includegraphics[width= \linewidth]{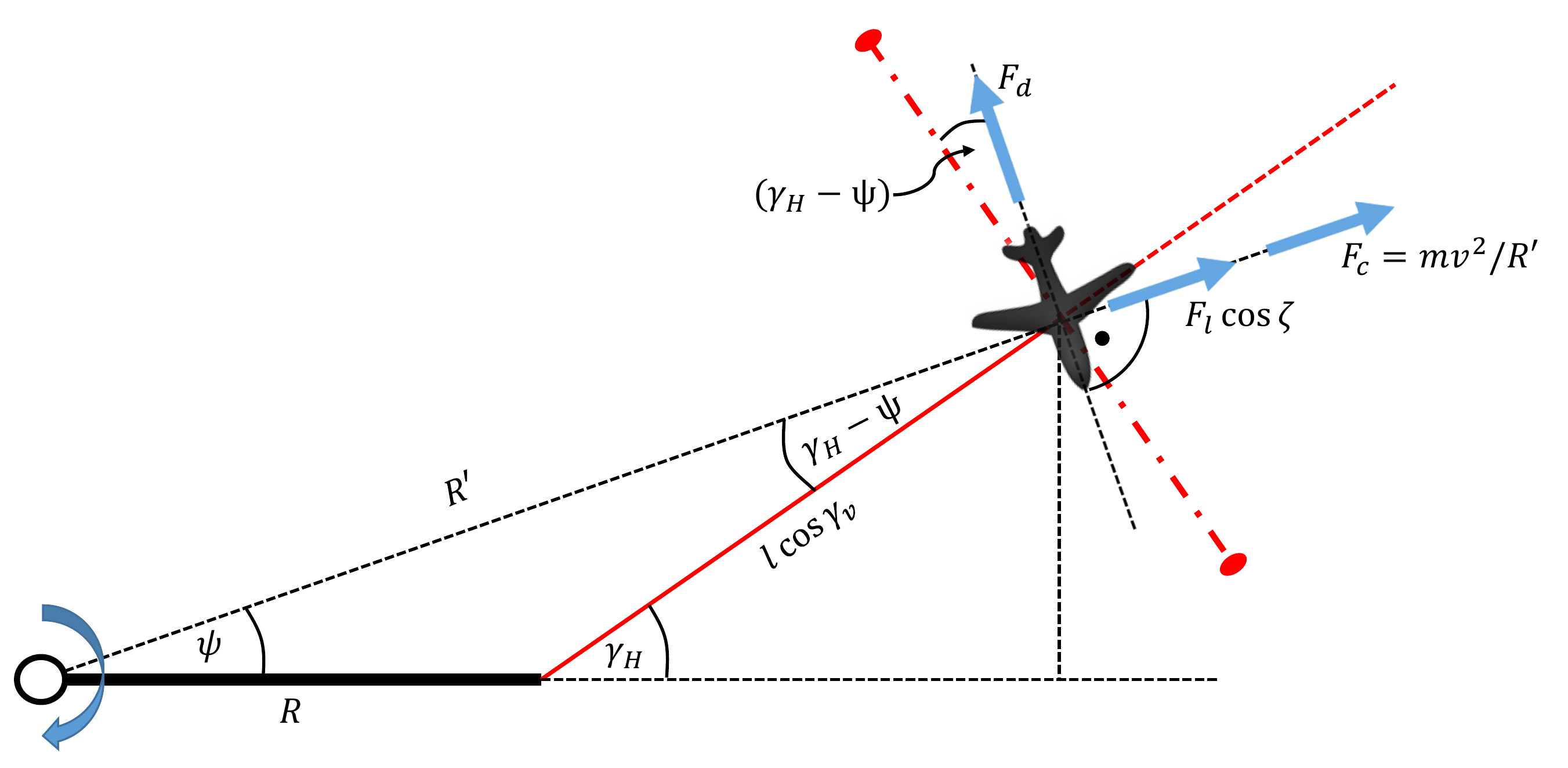}
		\caption{Drag, lift and centrifugal forces (or their components, respectively) and angles during the rotational take-off in the plane of the rotating arm. The rotating arm has a length $R$ and the tether (in red) of $l$.\label{Fig_RotStart1}}
	\end{center}
\end{figure}

\begin{figure}[t]
	\begin{center}
		\includegraphics[width=\linewidth]{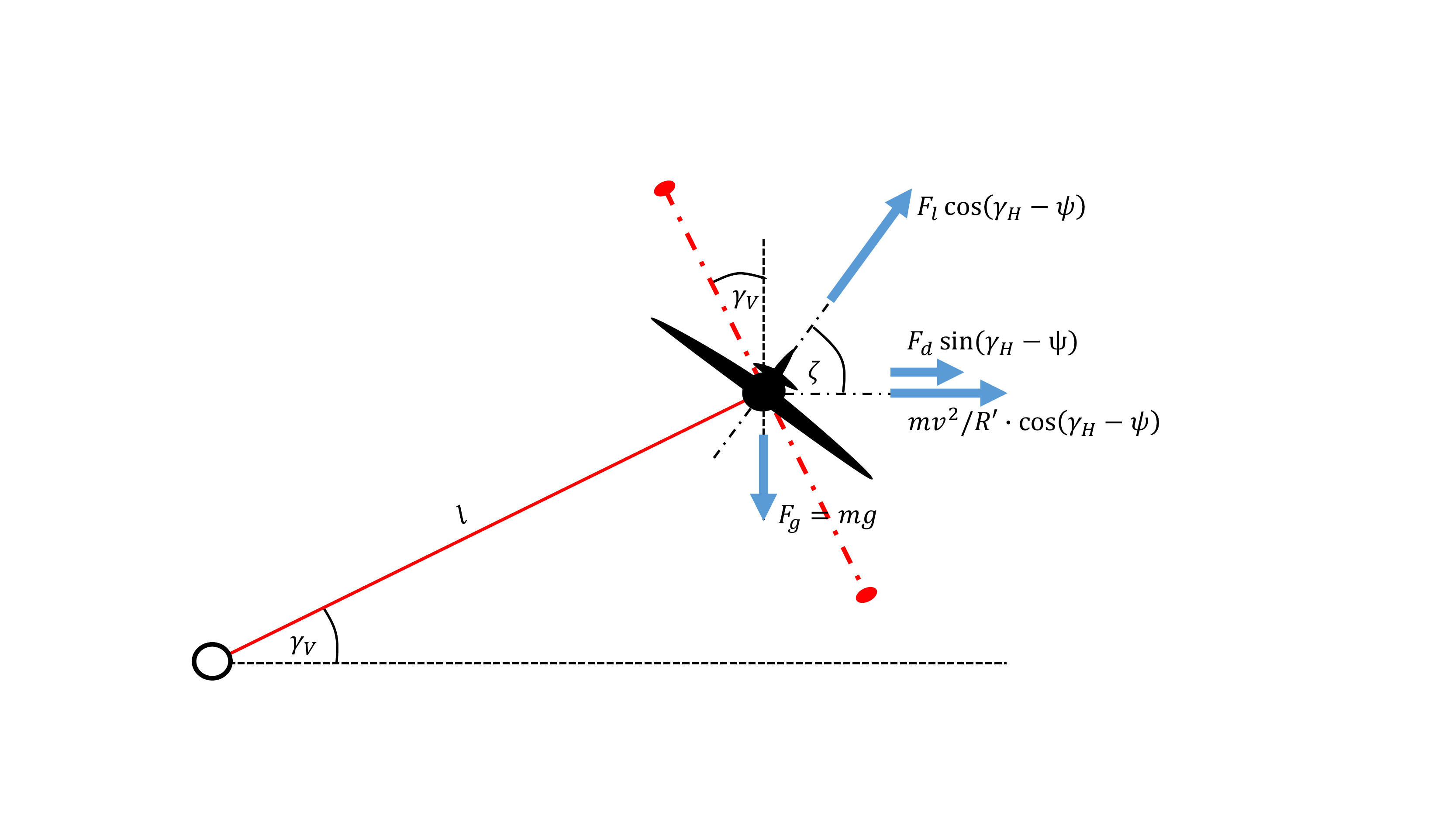}
		\caption{Drag, lift, centrifugal and gravitation forces (or their components, respectively) and angles during the rotational take-off in the plane perpendicular to the rotating arm and containing the tether.\label{Fig_RotStart2}}
	\end{center}
\end{figure}

A schematic arrangement of the rotational take-off is shown in Figure~\ref{Fig_RotStart}: the hull of the aircraft is attached to the tip of a rotating arm with length $R$ via the tether. The two angles $\gamma_v$ and $\gamma_h$ describe the orientation of the tether, assumed straight, with respect to the arm. The combination of lift force and centrifugal force due to the rotation leads to a reel out of the tether and the rise of the plane. If we assume that the angles $\gamma_v$ and $\gamma_h$ are constant during the rotational take-off, the sum of all forces perpendicular to the tether must cancel each other. Then, the required power to rotate the whole system (neglecting the drag of the rotating arm) is
\begin{equation} \label{rot_power}
	P_{g,2} = R T_\bot \omega,
\end{equation}
where
\begin{equation} \label{E:rot_tension_perp}
T_\bot = T\cdot\sin{(\gamma_H)}\cos{(\gamma_V)}
\end{equation}
is the tether tension $T$ projected onto the plane of the rotating arm and perpendicular to it and $\omega$ is the angular velocity of the system.

We consider a projection of Figure~\ref{Fig_RotStart} onto the plane of the rotating arm, as depicted in Figure~\ref{Fig_RotStart1}. Given $R,\,\omega,\,\gamma_H,\,\gamma_V$ and line length $l$, we define the angle $\psi$ and the distance $R'$ as:
\begin{equation}\label{E:psi}
	\psi \doteq \arctan\left(\frac{l\cdot\cos{(\gamma_V)}\cdot \sin{(\gamma_H)}}{R + l\cdot\cos{(\gamma_V)}\cdot\cos{(\gamma_H)}}\right),
\end{equation}
\begin{equation}\label{E:R'}
R'\doteq\dfrac{R+l\cos(\gamma_H)\cos(\gamma_V)}{\cos(\psi)}.
\end{equation}
Then, assuming that the absolute wind speed is zero, the aircraft will develop a lift force $F_l$ and a drag force $F_d$ whose magnitudes are equal to
\begin{equation} \label{E:lift_drag}
\begin{array}{rcl}
F_l&=&\dfrac{1}{2}\rho\,A\,C_l(R'\,\omega)^2\\
F_d&=&\dfrac{1}{2}\rho\,A\,C_{d,eq}(R'\,\omega)^2\\
\end{array}
\end{equation}
Figure~\ref{Fig_RotStart1} also shows the projections of all the considered forces (lift, drag, and centrifugal force) onto the plane of the rotating arm. The components perpendicular to the tether are the ones parallel to the dot-dashed line in the Figure. The requirement that they cancel each other yields
\begin{equation}\label{cond1}
	F_d \cos\left(\gamma_H-\psi \right) = \left(F_l \cos{(\zeta)} + m\frac{v^2}{R'}\right) \cdot \sin\left(\gamma_H-\psi\right),
\end{equation}
where $\zeta$ is the roll angle of the aircraft, as shown in Figure~\ref{Fig_RotStart2} which is the projection of Figure~\ref{Fig_RotStart} onto the plane perpendicular to that of the rotating arm and containing the tether. Again, the forces perpendicular to the tether are the ones parallel to the dot-dashed line in Figure~\ref{Fig_RotStart2}. Thus, the following condition must hold at the equilibrium, too:
\begin{equation}\label{cond2}
	\begin{split}
		&F_l \cos \left(\gamma_H-\psi \right) \sin\left(\zeta-\gamma_V \right) = mg\cdot \cos{(\gamma_V)} \\
		&+\left( m\frac{v^2}{R'}\cos\left(\gamma_H-\psi\right)+F_d \sin\left(\gamma_H-\psi\right)\right)\cdot \sin{(\gamma_V)}.
	\end{split}
\end{equation}
Finally, the tether tension in Eq.~\eqref{rot_power} is
\begin{equation}\label{E:Rot_tension}
	\begin{split}
		T = &F_l \cdot\cos\left(\gamma_H-\psi\right)\cos\left(\zeta-\gamma_V\right)-mg\cdot\sin{(\gamma_V)} \\
		&+ \left[F_d\sin\left(\gamma_H-\psi\right)+m\frac{v^2}{R'}\cos\left(\gamma_H-\psi\right)\right]\cdot\cos{(\gamma_V)}
	\end{split}
\end{equation}

Eqs.~\eqref{rot_power}-\eqref{E:Rot_tension} can be used to derive the power and ground area required for the rotational take-off. Since there exist many potential solutions that satisfy the equilibrium constraints \eqref{cond1}-\eqref{cond2}, we choose to evaluate this take-off approach by means of numerical optimization. We compute the involved variables (i.e. $\omega,\,\zeta$ etc.) and minimize the required mechanical power installed on the ground, $\overline{P}_{g,2}$, under certain operational constraints. More specifically, we fix the value of the arm length $R$ and, for each pair $(l,\gamma_V)$, we solve the following nonlinear program:
\begin{subequations}\label{E:rot_opt_prob}
\begin{align}
P^*_{g,2}(l,\gamma_V,R)&=&\min\limits_{\zeta,\omega,\gamma_H} (R T_\bot \omega)\label{E:rot_opt_cost}\\
&&\text{subject to}\nonumber\\
&&\textrm{Eqs. } \eqref{E:rot_tension_perp}-\eqref{E:Rot_tension}\label{E:rot_constr-overall}\\
&&\textrm{and } |\zeta-\gamma_V|\leq\overline{\zeta}\label{E:zeta_constr}
\end{align}
\end{subequations}
where the constraint \eqref{E:zeta_constr} is used to guarantee that the roll angle of the aircraft is such that the inner wing does not get too close to the tether, which might lead to entanglement and subsequent crash. Then, for each considered arm length $R$, we compute the peak required power as
\begin{equation}\label{E:peak_rot_R}
\overline{P}^*_{g,2}(R)=\min\limits_{\gamma_V\in[\underline{\gamma}_V,\,\overline{\gamma}_V]}\;\max\limits_{l\in[0,\,\overline{l}]}P^*_{g,2}(l,\gamma_V,R).
\end{equation}
The intervals $[\underline{\gamma}_V,\,\overline{\gamma}_V]$ and $[0,\,\overline{l}]$ considered in Eq.~\eqref{E:peak_rot_R} cover the range of reasonable equilibrium configurations that can occur when setting a constant vertical inclination $\gamma_V$ and reeling out the line. In particular, we assume that the line is reeled-out at a constant speed $v_l\ll \omega R'$, and that a specified vertical velocity $v_{c}$ of the aircraft is achieved. Then, from geometrical considerations we have that a minimum angle $\underline{\gamma}_V=\arcsin\left(\dfrac{v_c}{v_l}\right)$ shall be achieved.

The rationale behind problems \eqref{E:rot_opt_prob}-\eqref{E:peak_rot_R} is the following: For a given arm length $R$, we fix the vertical inclination of the line $\gamma_V$ during the ascend and we compute the required peak power over a reasonable range of line length values. Then, we search for the vertical inclination that achieves the lowest peak power. In this way, we obtain the minimal peak power, $\overline{P}^*_{g,2}(R)$, achievable with the considered arm length $R$ and the strategy of ascending with constant vertical inclination. Finally, we repeat this procedure over a range of arm lengths $R\in[\underline{R},\,\overline{R}]$ in order to find the minimal peak power $\overline{P}_{g,2}$ required to compute our quantitative criterium \textbf{C1}:
\begin{equation}\label{E:peak_rot_P}
\overline{P}_{g,2}=\min\limits_{R\in[\underline{R},\,\overline{R}]}\overline{P}^*_{g,2}(R).
\end{equation}
Regarding the required peak onboard power $\overline{P}_{ob,2}$ and additional mass $\Delta m$, both these quantities are virtually zero in this approach. Finally, the required ground area $A_{g,2}$ is equal to $\pi\,R_\text{opt}^2$, where $R_\text{opt}$ is the argument that minimizes \eqref{E:peak_rot_P}.


\subsection{Linear take-off with on-board propellers}\label{SS:linear}

In the following discussion of the linear take-off, we first analyze the on-ground acceleration phase and then the climbing phase.

\subsubsection{Acceleration phase on the ground}\label{linear}

The acceleration phase on the ground lasts until the take-off speed $v^*$ is reached:
\begin{equation}\label{E:takeoff_speed}
v^*=\sqrt{\dfrac{2(m+\Delta m_3) g}{\rho A C_l}},
\end{equation}
computed by setting $F_l=(m+\Delta m_3)\,g$ and using $F_l=\dfrac{1}{2}\rho A C_l v^{*2}$. Assuming that this speed shall be reached after a horizontal acceleration distance $L$, the required acceleration is $a={v^*}^2/(2L)$. The corresponding required force is then $F_\mathrm{g}=(m+\Delta m_3)\, a$. The other forces acting at take-off are significantly smaller, but not negligible, namely the drag force $F_d =\dfrac{1}{2} \rho C_{d,eq} A {v^*}^2$ and the viscous resistance $F_\mathrm{v} = c_\mathrm{v}\, v^*$, where $c_\mathrm{v}$ is the viscous friction coefficient of the system employed for the linear acceleration. Hence, the required maximal power on the ground is
\begin{equation}\label{E:takeoff_power}
	\overline{P}_{g,3} = v^{*}\,\left(F_\mathrm{g}+F_d+F_\mathrm{v}\right).
\end{equation}
As regards the land occupation, we choose to fix the travel length, such that it is independent from the wing size, and we assume that the system shall be able to adapt to the widest possible range of prevalent wind conditions, i.e. the linear acceleration phase can be carried out in all directions. At the same time, like the vertical take-off the area spanned by the wings throughout the ground launching phase is considered to be occupied by the system. Thus, we obtain \begin{equation}\label{E:linear_area}
A_{g,3}\simeq\dfrac{\pi L^2}{4}+\dfrac{\pi \lambda}{4}A.
\end{equation}

\subsubsection{Powering the plane during the ascend}\label{thrust_ascend}

\begin{figure}[t]
	\begin{center}
		\includegraphics[width=8cm]{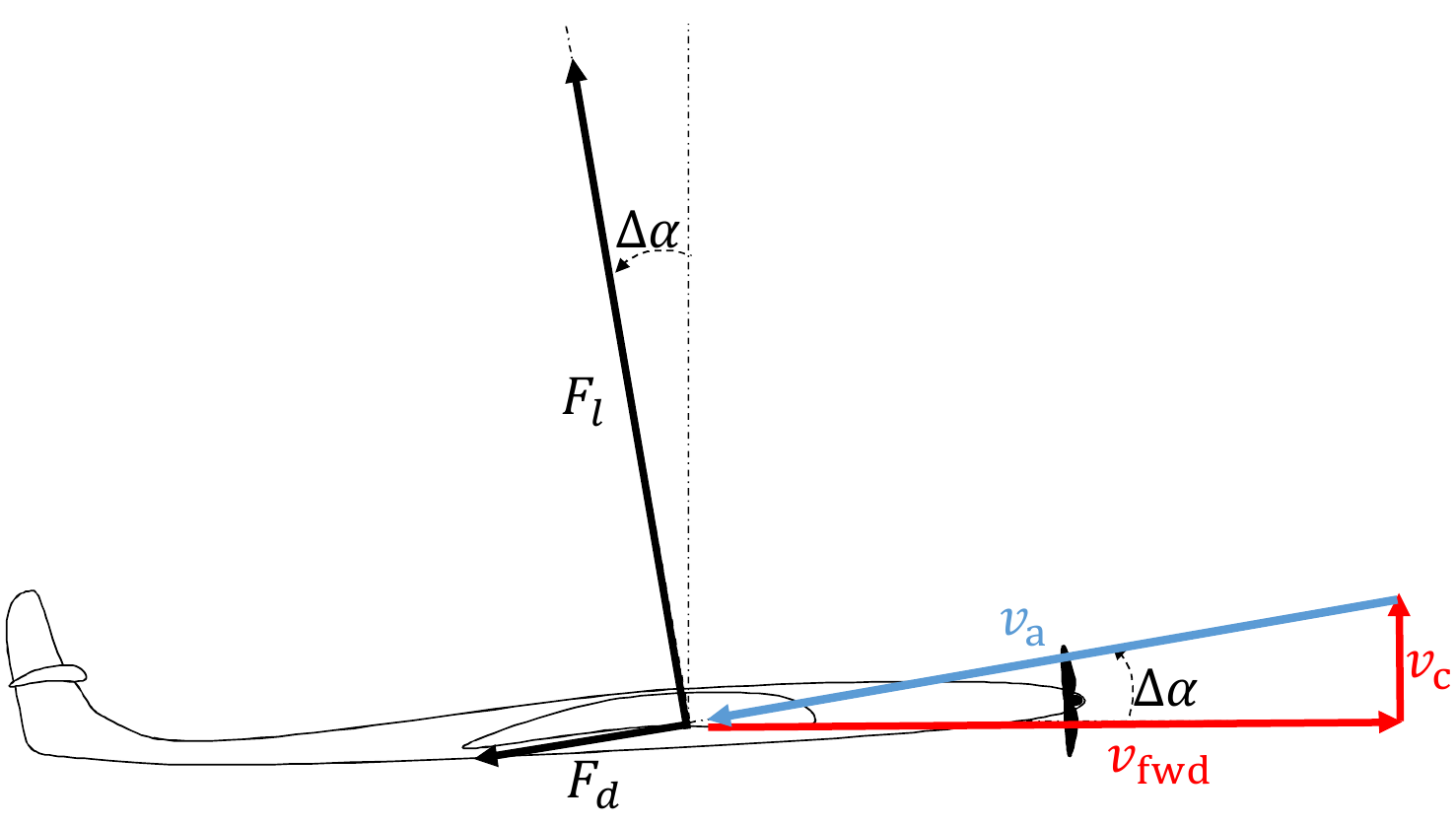}
		\caption{Schematic representation of an airplane with horizontal speed of $v_{\text{fwd}}$ (assuming no wind) and a vertical speed of $v_{\text{c}}$. The lift force has a component opposite to the thrust and the drag force has a component which adds to the gravitational pull. \label{Fig_LinStart}}
	\end{center}
\end{figure}

It is rather complicated for an external device like the winch to power the plane during the ascend. Indeed, on-board propellers are preferable because they can be small, since they do not have to accelerate the plane any further, and they shall just balance the aerodynamic drag and part of the lift depending on the climbing angle. In the following, we analyze the climbing phase assuming the worst conditions possible, i.e. with zero prevalent wind speed, which yields the peak on-board power. In the presence of wind, the climb may be carried out upwind with correspondingly lower power.

We denote  the vertical climb velocity with $v_\mathrm{c}$ again, see Figure~\ref{Fig_LinStart}. At the same time, the airplane moves horizontally with the speed $v_\mathrm{fwd}$ so that the total speed relative to the air is $v_\mathrm{a}=v_\mathrm{fwd}\cdot \sqrt{1+c_\mathrm{r}^2}$ with the climb ratio $c_\mathrm{r}:=v_\mathrm{c}/v_\mathrm{fwd}$. From Figure~\ref{Fig_LinStart}, it follows that $\sin{(\Delta\alpha)}=c_\mathrm{r}/\sqrt{1+c_\mathrm{r}^2}$ and $\cos{(\Delta\alpha)}=1/\sqrt{1+c_\mathrm{r}^2}$.

The vertical component of the lift force must counteract the gravitational pull and the vertical component of the drag force in order to yield a constant climb rate; i.e. the vertical equilibrium condition is $F_l\cdot \cos{(\Delta\alpha)}-F_d\cdot \sin{(\Delta\alpha)} = (m+\Delta m_3)g$. This gives
\begin{equation}\label{vert_equib}
	\frac{1}{2}\rho A C_l \sqrt{1+c_\mathrm{r}^2} \left(1-c_\mathrm{r} \frac{C_{d,eq}}{C_l} \right) v_\mathrm{fwd}^2 = (m+\Delta m_3) g.
\end{equation}	

For horizontal propulsion the required thrust is equal to the sum of the horizontal components of the lift and drag force, i.e.
\begin{equation}\label{hor_thrust}
	\begin{split}
		F_T &= F_l \cdot \sin{(\Delta\alpha)} + F_d \cdot \cos{(\Delta\alpha)} \\
		&= \frac{1}{2}\rho A C_l \sqrt{1+c_\mathrm{r}^2} \left(c_\mathrm{r} + \frac{C_{d,eq}}{C_l} \right) v_\mathrm{fwd}^2.
	\end{split}
\end{equation}

Considering that the climb ratio is typically of the order of 0.1-0.2 and that the aerodynamic efficiency of the aircraft is of the order of 10-20, we assume that $C_l/C_{d,eq} \gg c_\mathrm{r}$ and obtain from Eqs.~\eqref{vert_equib} and \eqref{hor_thrust} the final expression for the required thrust:
\begin{equation} \label{req_thrust}
	F_T = (m+\Delta m_3) g\cdot \frac{1+c_\mathrm{r}\frac{C_l}{C_{d,eq}}}{\frac{C_l}{C_{d,eq}}-c_\mathrm{r}} \approx (m+\Delta m_3) g\cdot \left(\frac{C_{d,eq}}{C_l}+c_\mathrm{r}\right).
\end{equation}

The required horizontal (forward) velocity can be calculated from \eqref{vert_equib}. Thus, for a desired climb rate $c_\mathrm{r}$, both force and horizontal velocity can be computed using Eqs.~\eqref{vert_equib} and \eqref{req_thrust}. Similarly to what discussed for the vertical take-off, the corresponding required peak power $\overline{P}_{ob,3}$ for the propellers is then given by:
\begin{equation}\label{lin_ob_power}
	\overline{P}_{ob,3} = \dfrac{F_T}{\eta}\left(\sqrt{\frac{F_T}{2\rho A_\mathrm{prop}} + \frac{v^2_\mathrm{fwd}}{4}} + \dfrac{1}{2} v_\mathrm{fwd}\right).
\end{equation}
For the propeller area $A_\mathrm{prop}$, we consider two propellers (this time with horizontal axis) with a diameter of half the chord and an efficiency of $\eta$.

Finally, as regards the additional on-board mass $\Delta m_3$, similarly to the vertical take-off we consider the energy density of on-board batteries and electric motors (see \eqref{E:heli_deltam}) and solve the resulting system of equations to obtain consistent values of $\overline{P}_{g,3},\,\overline{P}_{ob,3}$ and $\Delta m_3$.

%
%
%

\subsection{Results and Discussion}\label{SS:conclusion}

In this section, we apply the results presented so far to evaluate the criteria \textbf{C1}-\textbf{C3}. In  particular, we consider three different wing sizes and corresponding design parameters as shown in Table~\ref{wing_sizes}. The obtained results are used, together with the qualitative criteria \textbf{C4}-\textbf{C5}, to discuss the considered take-off approaches and draw conclusions on their viability. For the computation of \textbf{C1}, the mechanical power $P^*_m$ is calculated with Eq.~\eqref{E:power2} with a wind speed $W=15\,$m/s. Regarding the energy density of on-board batteries and the power density of on-board motors, we considered $E_\text{batt}=720\,\textrm{kJ/kg}$ and $E_\text{mot}=2.5\,$kg/kW \cite{Bont10}. Table~\ref{wing_sizes} also shows in bold the results obtained according to the analysis described in sections \ref{SS:rotorcraft}-\ref{SS:linear}, including the values of $\overline{P}_{g,i},\,\overline{P}_{ob,i},\,\Delta m_i,$ and $A_{g,i}$. Finally, Table~\ref{T:results_power_mass_area} summarizes the values of the scaling factors that define the criteria \textbf{C1}-\textbf{C3}, obtained with the parameters of Table~\ref{wing_sizes}. Before drawing a final assessment, we briefly comment on the results obtained with each approach.

\begin{table*}\center
	\caption{Design parameters for the assessment of the different take-off concepts. Bold-faced parameters are the results obtained according to the assumptions and  analysis described in sections \ref{SS:rotorcraft}-\ref{SS:linear}.\label{wing_sizes}}
	\begin{tabular}{|l|c|c|c|}
		\hline
        \textbf{Parameter}& \textbf{Aircraft 1} & \textbf{Aircraft 2} & \textbf{Aircraft 3} \\ \hline\hline
        \multicolumn{4}{|c|}{Common parameters}\\\hline\hline
				Wing span $d$ (m) & 5 & 10 & 20 \\ \hline
        Aspect ratio $\lambda$ & \multicolumn{3}{|c|}{10}\\ \hline
				Chord $d/\lambda$ (m) & 0.5 & 1 & 2 \\\hline
				Wing area $A$ (m$^2$) & 2.5 & 10 & 40 \\ \hline
        Wing loading $w_{l} = m/A$ (kg/m$^2$)&  \multicolumn{3}{|c|}{15}\\ \hline
				Mass $m_0$ (without additional equipment for take-off & 37.5 & 150 & 600 \\ \hline
				Lift coefficient $C_l$ & \multicolumn{3}{|c|}{1}\\ \hline
				Drag coefficient $C_{d,eq}$ & \multicolumn{3}{|c|}{0.1}\\ \hline
				Desired vertical velocity $v_{\textrm{c}}$ (m/s) & \multicolumn{3}{|c|}{1} \\ \hline
				Propeller efficiency $\eta$ & \multicolumn{3}{|c|}{0.7} \\ \hline
				Peak mechanical power $P^*_m$ with  $W=15$ m/s (kW) & 75 & 300 &1200 \\ \hline\hline
        \multicolumn{4}{|c|}{Vertical take-off}\\\hline\hline
        Target height $h$ (m) & \multicolumn{3}{|c|}{100} \\ \hline
        Energy density of on-board batteries $E_\text{batt}$ (kJ/kg) & \multicolumn{3}{|c|}{720}\\ \hline
        Power density of on-board motors $E_\text{mot}$ (kW/kg) &\multicolumn{3}{|c|}{2.5}\\\hline
        Propeller diameter $d/\lambda$ (m) & 0.5 & 1 & 2 \\\hline
        Peak additional on-board power $\overline{P}_{ob,1}$ (kW) & \textbf{14} & \textbf{56}& \textbf{223} \\\hline
        Additional on-board mass $\Delta m_1$ (kg) & \textbf{8} & \textbf{30} & \textbf{120} \\\hline
        Required ground area $A_{g,1}$ (m$^2$) & \textbf{20} & \textbf{80} & \textbf{315} \\
        \hline\hline
        \multicolumn{4}{|c|}{Rotational take-off}\\\hline\hline
        Maximum angle between the wings & \multicolumn{3}{|c|}{ }\\ and the plane perpendicular to the line $\overline{\zeta}$ (deg) &  \multicolumn{3}{|c|}{50}\\\hline
        Reel-out speed of the line $v_l$ (m/s)&  \multicolumn{3}{|c|}{1.6} \\\hline
        Minimum vertical inclination $\underline{\gamma}_V$ (deg) &  \multicolumn{3}{|c|}{\textbf{40}}\\\hline
        Maximum vertical inclination $\overline{\gamma}_V$ (deg) &  \multicolumn{3}{|c|}{90}\\\hline
        Minimum arm length $\underline{R}$ (m) & \multicolumn{3}{|c|}{\textbf{30}}\\\hline
        Maximum arm length $\overline{R}$ (m) & \multicolumn{3}{|c|}{50}\\\hline
        Optimal arm length $R_\text{opt}$ (m) & \multicolumn{3}{|c|}{\textbf{50}} \\ \hline
        Maximal angular velocity $\omega$ (rad/s)&  \multicolumn{3}{|c|}{\textbf{0.4}}\\\hline
        Maximal tangential velocity of the tip of the arm $\omega\,R$ (m/s) &  \multicolumn{3}{|c|}{\textbf{20}}\\\hline
        Peak additional ground power $\overline{P}_{g,2}$ (kW) & \textbf{3} & \textbf{12} & \textbf{47} \\\hline
        Additional on-board mass $\Delta m_2$ (kg) & \multicolumn{3}{|c|}{0} \\\hline
        Required ground area $A_{g,2}$ (m$^2$) & \multicolumn{3}{|c|}{\textbf{7854}}\\\hline\hline
        \multicolumn{4}{|c|}{Linear take-off with on-board propellers}\\\hline\hline
        Ground travel distance $L$ (m) & \multicolumn{3}{|c|}{12} \\ \hline
				Target height $h$ (m) & \multicolumn{3}{|c|}{100} \\ \hline
        Viscous friction coefficient $c_v$ (kg/s) & 0.1 & 0.3 & 1 \\ \hline
        Take-off speed $v^*$ (m/s) & \multicolumn{3}{|c|}{\textbf{15.7}} \\ \hline
        Propeller's diameter $d/(2\lambda)$ (m) & 0.25 & 0.5 & 1 \\ \hline
        Peak additional ground power $\overline{P}_{g,3}$ (kW) & \textbf{8} & \textbf{31} & \textbf{124} \\\hline
        Peak additional on-board power $\overline{P}_{ob,3}$ (kW) &\textbf{2} & \textbf{9} & \textbf{37} \\\hline
        Additional on-board mass $\Delta m_3$ (kg) & \textbf{2} & \textbf{5} & \textbf{20} \\\hline
        Required ground area $A_{g,3}$ (m$^2$) & \textbf{132} & \textbf{192} & \textbf{428} \\\hline
	\end{tabular}
\end{table*}

{\bf Vertical take-off.} As expected, this approach requires the largest amount of additional on-board power (about 20\% of the peak mechanical power of the system) and of additional mass (20\% of the aircraft mass), see Table \ref{T:results_power_mass_area}. On the other hand, the required ground area turns out to be the smallest among the three approaches. The additional complexity (criterium \textbf{C4}) can be substantial, since the aircraft and on-board equipment have to be designed to sustain the large accelerations experienced during crosswind flight, and since large electric on-board power is required. This might require a completely new design of the wing. The additional mass also leads to a larger cut-in speed for the generator, since a larger wind speed will be required for the system to be able to remain airborne during power generation. Moreover, in a deeper analysis the presence of the propellers will have a detrimental influence on the aerodynamics, hence either requiring a larger wing for the same power, or giving lower power for the same size. These aspects lead in turn to a reduced capacity factor. The possibility to take-off in a large range of wind conditions (criterium \textbf{C5}) is in principle given, although more detailed studies should be carried out to assess whether the control surfaces and the propellers can effectively stabilize the aircraft during the ascend with relatively strong wind. 

\begin{table*}\center
	\caption{Results for the quantitative performance criteria \textbf{C1}, \textbf{C2}, and \textbf{C3} (Eqs.\eqref{E:crit_power}-\eqref{E:ground_area}~) with the parameters of Table \ref{wing_sizes}. \label{T:results_power_mass_area}}
	\begin{tabular}{|l|l|c|c|c|c|c|}
		\hline
		\multicolumn{2}{|c|}{}& \multicolumn{2}{|c|}{\textbf{C1}: power} & \multicolumn{1}{|c|}{\textbf{C2}: mass} & \multicolumn{2}{|c|}{\textbf{C3}: area} \\ \hline
		\multicolumn{2}{|l|}{Concept}   & $\eta_{P_{g},i}$ (\%) &$\eta_{P_{ob},i}$ (\%)  &  $\eta_{m,i}$ (\%)   & $\underline{A}_{g,i}$ &  $\eta_{A_g,i}$ (\%) \\ \hline\hline
		\multicolumn{2}{|l|}{Vertical} & 0  & 19  & 21  & 0 &  $\frac{\pi\lambda}{4}$ \\ \hline\hline
		\multicolumn{2}{|l|}{Rotational} &4  & 0 &0 & $\frac{\pi R^2}{4}$& 0\\ \hline\hline
\multicolumn{2}{|l|}{Linear} & 11  & 3  & 5  & $\frac{\pi L^2}{4}$& $\frac{\pi\lambda}{4}$ \\ \hline
 	\end{tabular}
\end{table*}

\begin{figure}[!htb]
 \begin{center}
  \includegraphics[width=8cm]{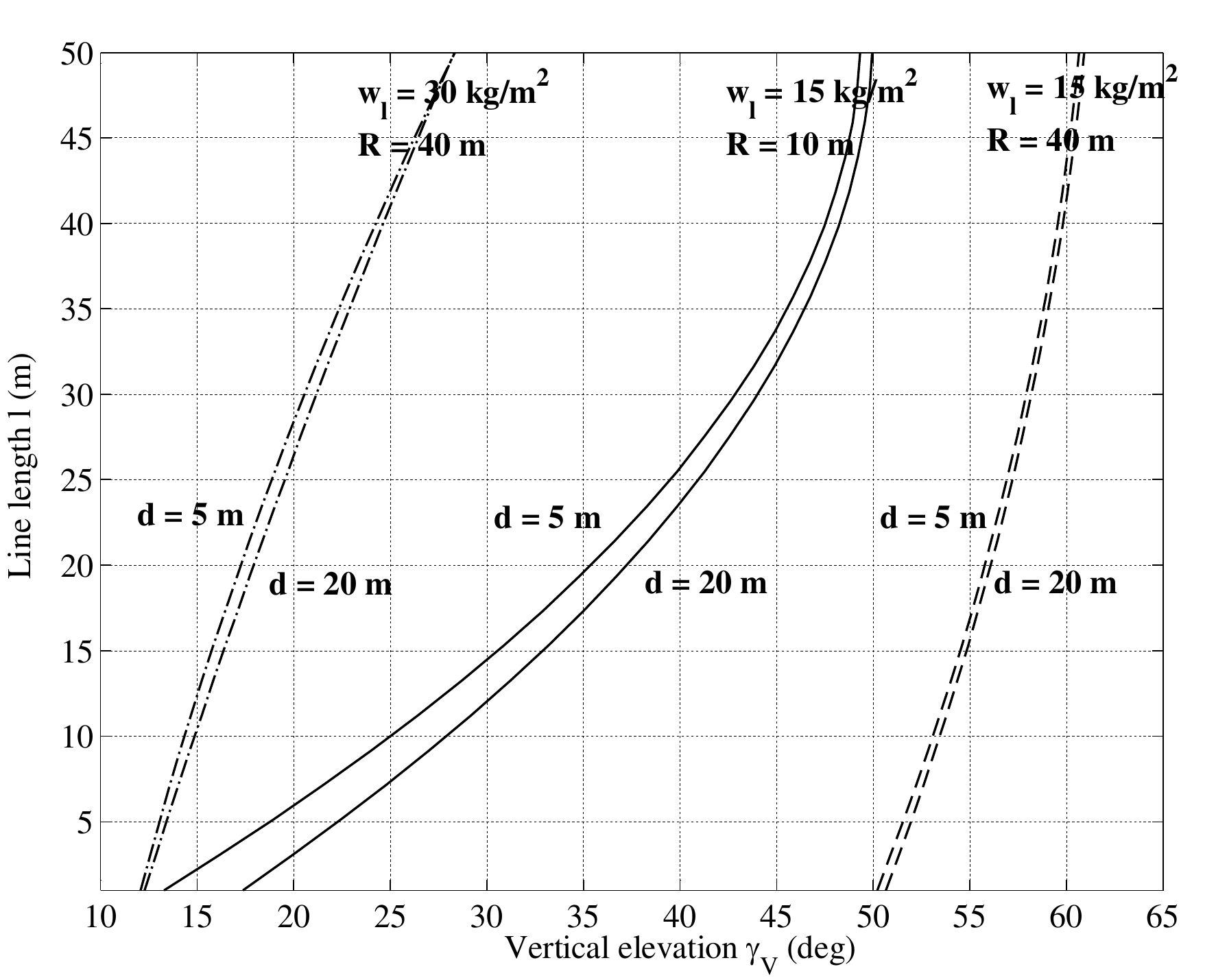}
  \caption{Analysis of the rotational take-off. Curves showing the maximum elevation angle $\gamma_V$ that can be achieved as a function of the line length with $w_l=15\,$kg/m$^2$ and $R=10\,\textrm{m}$ (solid lines), $w_l=15\,$kg/m$^2$ and $R=40\,\textrm{m}$ (dashed), and $w_l=30\,$kg/m$^2$ and $R=40\,\textrm{m}$ (dash-dotted). For each combination of $w_l$ and $R$, two values of wingspan ($d=5\,\textrm{m}$ and $d=20\,\textrm{m}$) are shown.}
  \label{F:minimum_R}
 \end{center}
\end{figure}

{\bf Rotational take-off.} While the results for the vertical and linear approaches are derived in a straightforward way from the equations presented in sections \ref{SS:rotorcraft} and \ref{SS:linear}, some more comments are due on the results pertaining to the rotational take-off. The application of the optimization procedure described in section \ref{SS:rotational} provides several interesting outcomes. First, it turns out that there exist a minimal arm length $\underline{R}$ that allows the system to achieve vertical inclination angles larger than the minimum required one, i.e. $\underline{\gamma}_V$. The value of $\underline{R}$ mainly depends on the wing loading $w_l$, while it is not affected significantly by the wing size, as shown in Figure \ref{F:minimum_R} which presents the curves of maximum $\gamma_V$ values that can be achieved as a function of line length $l$, for various combinations of wing loading $w_l$, arm length $R$ and wingspan $d$.

The main explanation for this phenomenon is that the aerodynamic forces have to counteract the centrifugal force (see section \ref{SS:rotational}), which decreases as the arm length $R$ increases. For the wing loading and minimal vertical inclination values chosen for our comparison, i.e. $w_l=15\,$kg/m$^2$ and $\gamma_V = 40^\circ$, we obtain $\underline{R}\simeq 30\,$m, as reported in Table \ref{wing_sizes}.

\begin{figure}[!htb]
 \begin{center}
  \includegraphics[width=8cm]{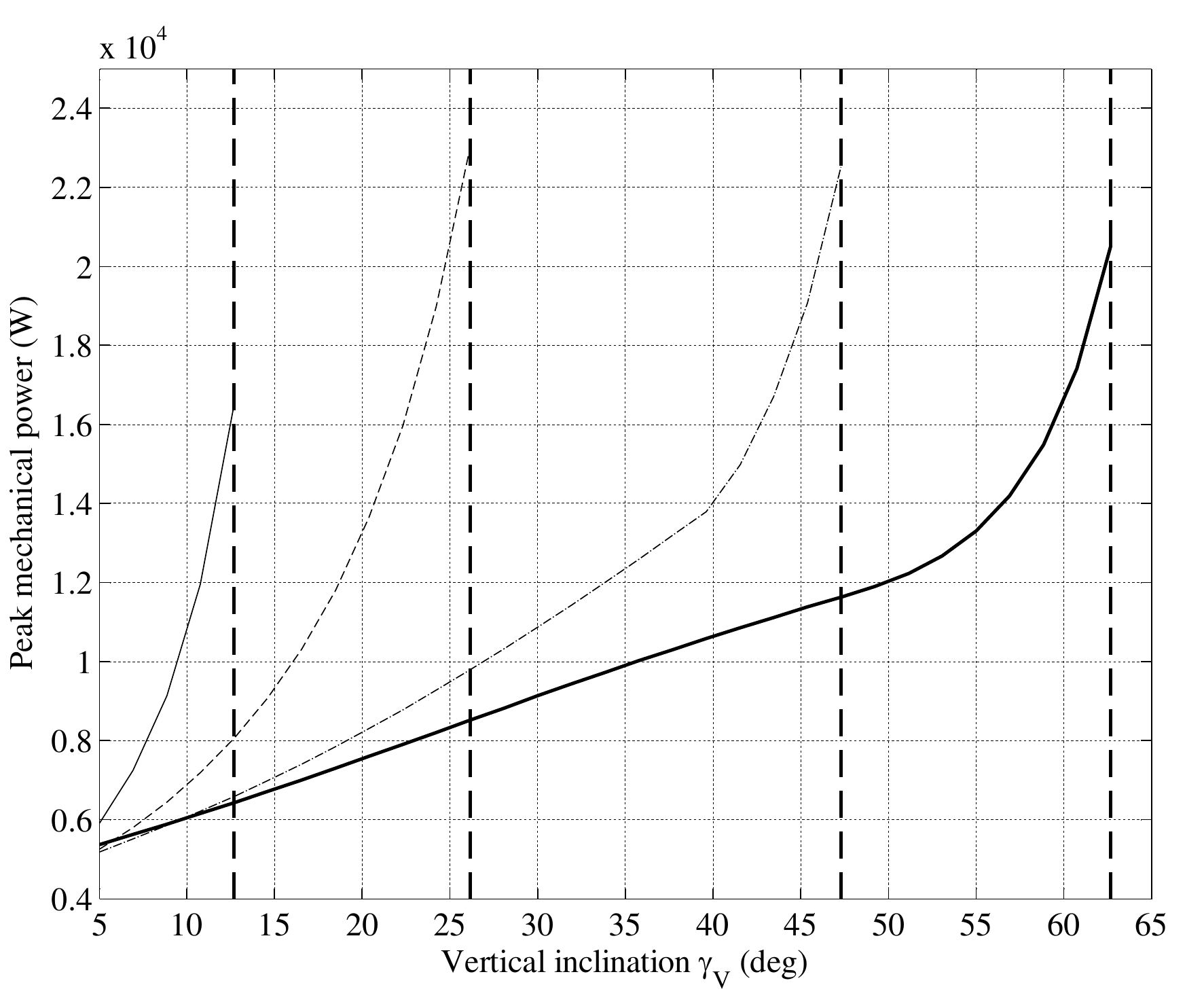}
  \caption{Analysis of the rotational take-off. Curves showing the peak ground mechanical power $P^*_{g,2}$ for $l=1\,$m as a function of the elevation angle $\gamma_V$ and with $R=10\,$m (thin solid line), $R=20\,$m (dashed), $R=40\,$m (dash-dotted) and $R=80\,$m (thick solid line). The vertical dashed lines indicate the maximum elevation angle achievable for each considered arm length. Wing span $d=10\,$m, wing loading $w_l=15\,$kg/m$^2$.}
  \label{F:rot_power_vs_gamm_v_result}
 \end{center}
\end{figure}

Second, the required peak power increases with $\gamma_V$, since the equilibrium conditions \eqref{cond1}-\eqref{cond2} become less favorable and a larger rotational speed is required to generate enough lift to maintain the desired vertical inclination. This is shown in Figure \ref{F:rot_power_vs_gamm_v_result}. Hence, for the sake of minimizing the required additional power, the minimum vertical inclination is chosen.

\begin{figure}[!htb]
 \begin{center}
  \includegraphics[width=8cm]{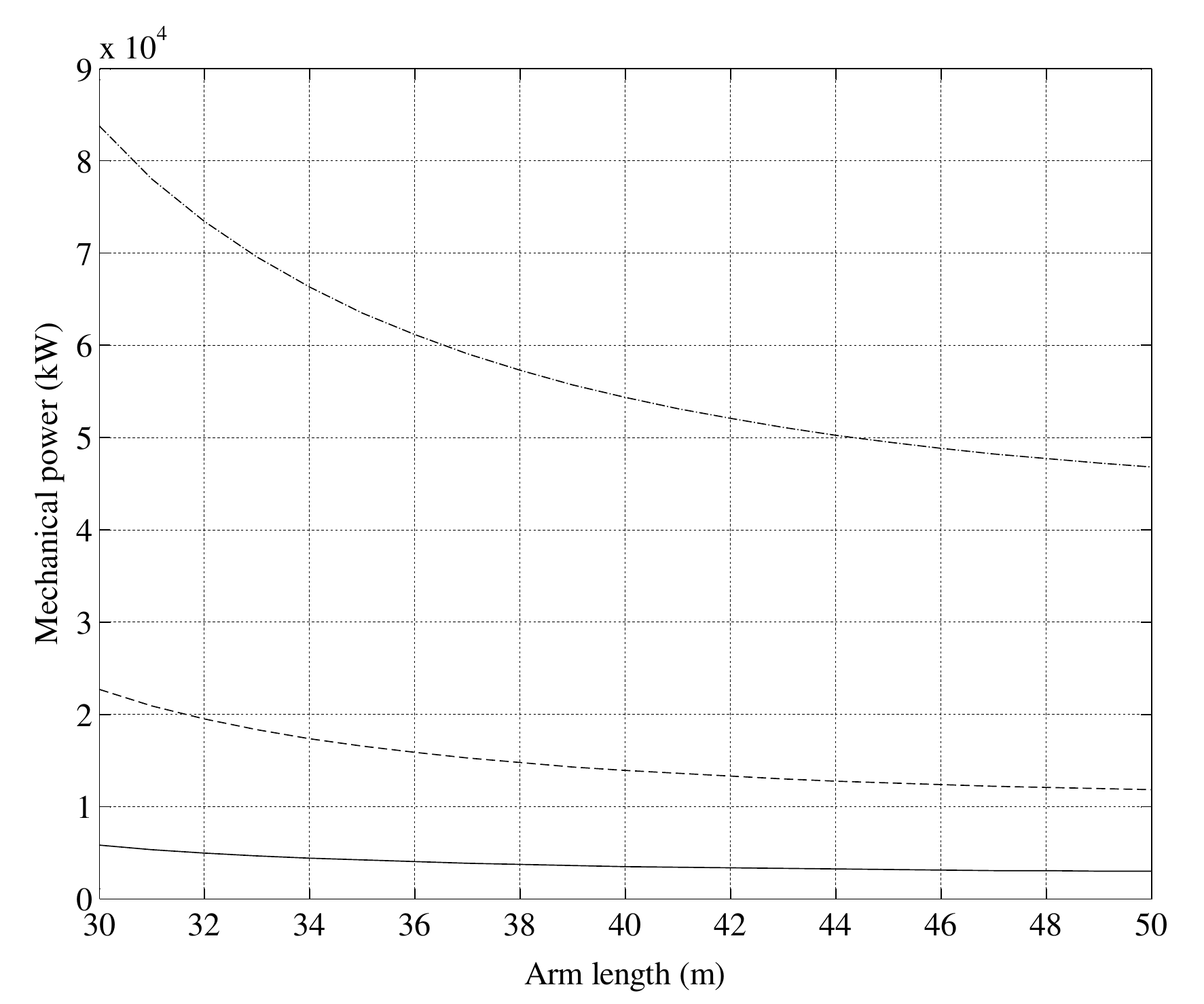}
  \caption{Analysis of the rotational take-off. Curves showing the peak ground mechanical power $P^*_{g,2}$ for $l=1\,$m as a function of the arm length $R$ for $d=5\,$m (solid), $d=10\,$m (dashed), and $d=20\,$m (dash-dotted). Wing loading $w_l=15\,$kg/m$^2$.}
  \label{F:rot_power_vs_R_result}
 \end{center}
\end{figure}

Last, the peak mechanical power decreases with the arm length and approaches an asymptotic value. The reason is that, as the centrifugal force decreases (i.e. $R$ increases), the aerodynamic forces have to win just the aircraft weight. This condition leads asymptotically,  for growing $R$, to a minimum required tangential speed and corresponding forces which then determine the required power to rotate the arm. The mentioned trend is shown in Figure \ref{F:rot_power_vs_R_result}. In order to restrict our analysis to a finite value of $R$, we chose an upper bound of $\overline{R}=50\,\textrm{m}$, which is then the optimal value according to Eq.~\eqref{E:peak_rot_P}.

Due to the mentioned findings, the rotational start-up results in the lowest value of peak ground and on-board power, as well as the lowest value of additional mass, but a very large ground-area occupation as compared with the other two approaches, see Tables \ref{wing_sizes}-\ref{T:results_power_mass_area}. Such a land occupation is fundamentally linked to the wing loading as discussed above, i.e. it is not possible to decrease the land occupation below a minimum threshold by increasing the installed power or decreasing the wing size. As a matter of fact, the minimum ground occupation is quite large for a reasonable wing loading.

As regards complexity (\textbf{C4}), this is expected to be large, considering that the system would feature a 50-m-long rigid arm whose tip rotates at about 20 m/s. Moreover, the whole main winch should rotate as well with many full revolutions while at the same time reeling the line, which poses a challenge for the winch mechanics and the electrical connections. The manufacturing and installation costs of such a structure could be comparable to those of a traditional wind turbine and appear to be prohibitive for the economic viability of the approach. Finally, about wind adaptation (\textbf{C5}) it is unclear how this concept would handle a strong prevalent wind during take-off, when the relative wind speed could changes e.g. by $\pm$10 m/s during a half turn, with the aircraft speed relative to ground of about 20$\,$m/s.

{\bf Linear take-off.} The required peak power installed on the ground for this approach is larger than that of the rotational take-off, however with a significantly smaller required area. Moreover, differently from the rotational take-off, in the linear take-off the ground area and required power can be easily traded off. As regards the on-board power and additional mass, they result to be about six times smaller than for the vertical take-off. The required ground occupation is comparable to the vertical take-off and dominated by the wing size when scaling up, hence it turns out to be quite favorable. About the complexity of the approach, this appears to be small, since in principle one could envision a solution where the winch used to generate power is also employed in the initial phase of the take-off, e.g. by means of a clutch to (dis-) engage a linear motion system to accelerate the aircraft. Similarly, the on-board propellers and batteries are necessary in any case to power the on-board control systems, hence the use of slightly larger and more powerful on-board motors does not appear to be critical. Moreover, the on-board propellers can also be used to re-charge the batteries to supply energy to the control system during long periods of power generation. Finally, since the whole setup can be turned, the take-off is independent of the current prevalent wind direction.

{\bf Discussion.} The results presented so far indicate that both the vertical and the rotational take-off require extensive modifications of the AWE system, which will have a strong influence on the design and require significant additional equipment. On the other hand, the  linear take-off approach will have less impact on the system design. If the main winch can also be used for the acceleration phase, the additionally required equipment is in fact reduced to a minimum. In terms of mechanical power, the linear take-off provides a good tradeoff between on-board and on-ground power. Moreover the additional on-board components like batteries and small propellers will have further applications, like powering the on-board electronics. Finally, the land occupation of the linear take-off is almost as small as that of the vertical one. For these reasons, we favor the linear take-off for rigid-wing AWE systems with ground-based electric generation. This approach will be analysed in more detail in the following section.

\section{Simulation of a linear take-off approach}\label{S:simulation}
In this section, we further study, by means of numerical simulations, the linear take-off combining ground motors and on-board propellers. We first introduce a dynamical model of the system, then we describe the control algorithms to carry out the take-off maneuver, finally we present the simulation results and compare them with the static equations derived in section \ref{SS:linear}.

\subsection{A dynamical model for linear take-off}\label{SS:model}
\begin{figure}[!htb]
 \begin{center}
  \includegraphics[trim= 0cm 0cm 0cm 0cm,width=8cm]{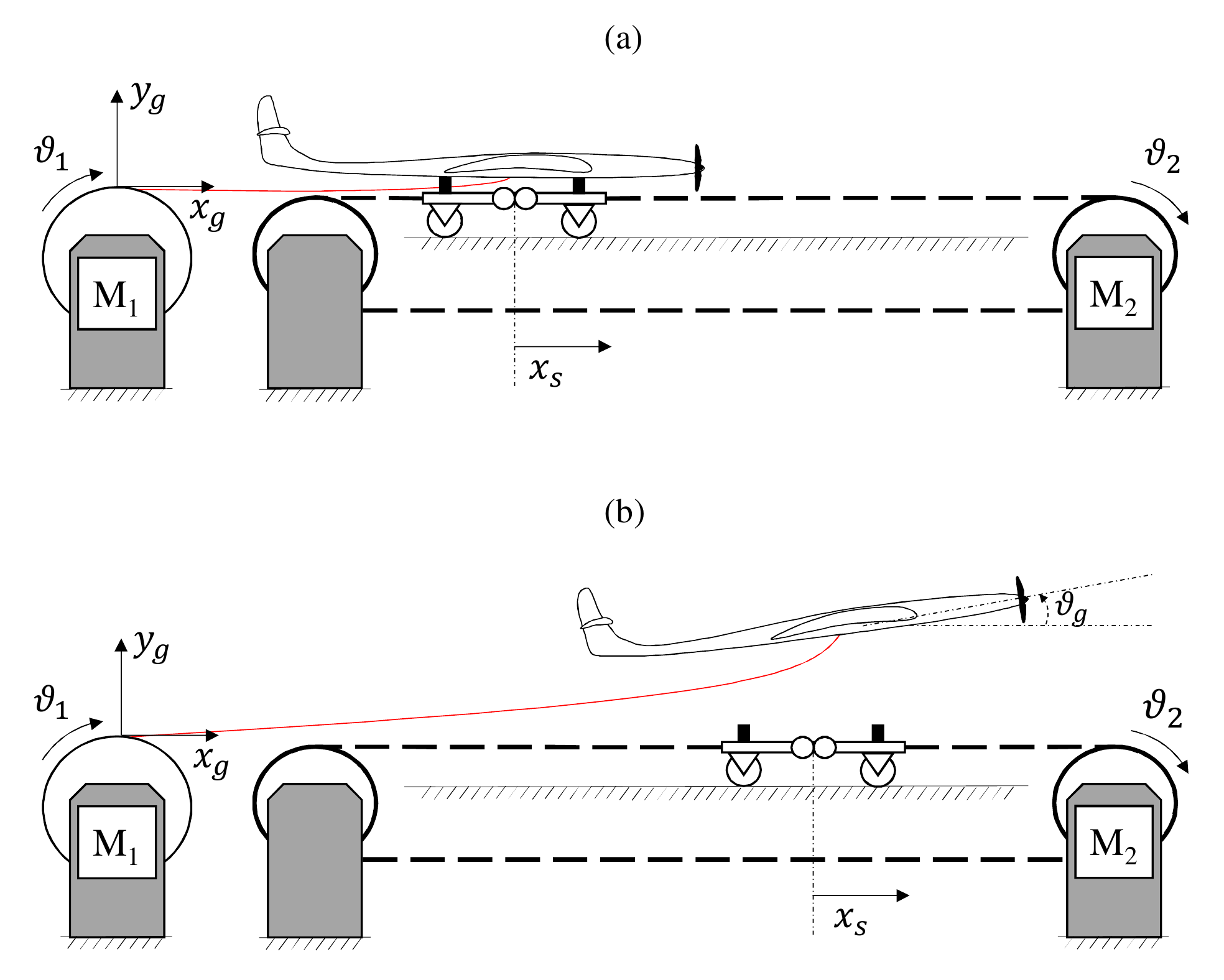}
  \caption{Sketch of the system considered to simulate the linear take-off procedure. (a) First operating mode, with the aircraft carried by the slide up to take-off speed; (b) second operating mode, with the aircraft gaining altitude by means of the on-board propeller.}
  \label{F:Dyn_system_sketch}
 \end{center}
\end{figure}

We consider a ground unit composed of a winch, where the aircraft's tether is coiled, and of a linear motion system, whose aim is to accelerate the aircraft up to take-off speed, see Figure \ref{F:Dyn_system_sketch}. The winch rotation is controlled by a  geared motor/generator $M_1$, which is the main electrical machine of the AWE system, responsible for converting mechanical power into electricity during the power generation cycles. The linear motion system consists of a slide, carrying the aircraft during take-off, that can move along rails. The slide motion is controlled by a second geared motor $M_2$ through a transmission system (e.g. a belt). The slide is equipped with sheaves that guide the tether from the winch to the attachment point on the aircraft. This system can be well described by a hybrid dynamical model: a first operating mode (Figure \ref{F:Dyn_system_sketch}(a)) describes the system's behavior from zero speed up to the take-off, when the aircraft and the slide can be considered as a unique rigid body; a second operating mode (Figure \ref{F:Dyn_system_sketch}(b)) describes the aircraft motion after take-off, when it is separated from the slide. For the sake of simplicity, we consider a two-dimensional motion only in the second mode, i.e. vertical and horizontal displacements and pitch rotation of the aircraft, assuming that suitable stabilizing systems act on the on-board actuators (rudder and ailerons) in order to keep the roll and yaw angles at small values, counteracting potential lateral wind turbulence. Moreover, we assume that no wind opposite to the take-off direction is present, i.e. the take-off is carried out only by means of the ground motors and on-board propellers. In case of substantial wind, we assume the system to be capable to orient the rails according to the wind direction to take advantage of the additional apparent wind velocity, hence reducing the take-off speed. Thus, the conditions simulated here provide the worst-case in terms of required power, in line with the analysis of section \ref{SS:linear}. All the equations presented in the following have been derived by applying Newton's second law of motion.

The state of the model, i.e. the variables that describe completely and univocally its configuration at any time instant $t$, is given by $x(t)\doteq[\vartheta_{M_1}(t),\,$ $\dot{\vartheta}_{M_1}(t),\,$ $\vartheta_{M_2}(t),\,$ $\dot{\vartheta}_{M_2}(t),\,$ $x_g(t),\,$ $\dot{x}_g(t),\,$ $y_g(t),\,\dot{y}_g(t),\,$ $\vartheta_{g}(t),\,\dot{\vartheta}_{g}(t)]^T$, where $\vartheta_{M_1}$ is the angular position of the winch, $\dot{\vartheta}_{M_1}\doteq\frac{d\vartheta_{M_1}}{dt}$ its angular speed, $\vartheta_{M_2},\,\dot{\vartheta}_{M_2}$ are the angular position and speed of the motor that controls the linear motion system,
$x_g(t),\,$ $\dot{x}_g(t),\,$ $y_g(t),\,\dot{y}_g(t),\,$ the horizontal ($x$) and vertical ($y$) positions and speeds of the aircraft's center of gravity in an inertial reference frame. The latter has its center at the point where the tether exits the winch,  the $x_g$-axis parallel to the ground and the $y_g$-axis vertical and pointing upwards (see Figure \ref{F:Dyn_system_sketch}). Finally, $\vartheta_{g}(t),\,\dot{\vartheta}_{g}(t)$ are the aircraft's pitch angle and its rate. The manipulated inputs available to control and operate the system are denoted with $u(t)\doteq[C_{M_1}(t),\,C_{M_2}(t),\,F_T(t)]^T$ where $C_{M_1},\,C_{M_2}$ are the torques applied by the two electrical machines, and $F_T$ is the thrust force exerted by the on-board propeller. The motor torques considered in the model are taken after any gear that can be installed between the motor and the winch (respectively the belt's pulley) to adapt the motor's torque/speed profile to the application. In the following, for the sake of simplicity we denote with $x_j$ (resp. $u_j$) the $j^\text{th}$ component of the state (resp. input) vector defined above. Assuming that the linear motion system is realized by a belt, driven by a pulley directly attached to the shaft of motor $M_2$, and neglecting its elasticity, 
the model is given by the following equations in the first operating mode:
\begin{equation}\label{E:dyn_mod_first}
\begin{array}{rcl}
\dot{x}_1(t)&=&x_2(t)\\
\dot{x}_2(t)&=&\dfrac{1}{J_{M_1}}(r_{M_1}\,T(t)-\beta_{M_1}\,x_2(t)+u_1(t))\\
\dot{x}_3(t)&=&x_4(t)\\
\dot{x}_4(t)&=&\dfrac{1}{J_{M_2}+(m_s+m)\,r_{M_2}^2}\,(r_{M_2}\,(-T(t)+\\
&
&-F_d(t)\,\cos(\Delta\alpha(t))+\\
& &+F_l(t)\,\sin(\Delta\alpha(t))-\beta_s\,r_{M_2}\,x_4(t))\\
& &-\beta_{M_2}\,x_4(t)+u_2(t))\\
\dot{x}_5(t)&=&x_6(t)\\
\dot{x}_6(t)&=&r_{M_2}\,\dot{x}_4(t)\\
\dot{x}_7(t)&=&x_8(t)\\
\dot{x}_8(t)&=&0\\
\dot{x}_9(t)&=&x_{10}(t)\\
\dot{x}_{10}(t)&=&0.\\
\end{array}
\end{equation}
In \eqref{E:dyn_mod_first}, $r_{M_1}$ is the radius of the winch (assuming for simplicity that the latter is directly connected to the motor/generator), $r_{M_2}$ the radius of the pulley that links motor $M_2$ to the belt, $J_{M_1},\,J_{M_2}$ the moments of inertia of the winch and of the pulley plus their respective motors, $\beta_{M_1},\,\beta_{M_2}$ their viscous friction coefficients, $m_s$ the mass of the slide, $\beta_s$ the viscous friction coefficient of the belt/slide/rail system, $m$ the mass of the aircraft. The angle $\Delta\alpha$ is defined as:
\begin{equation}\label{E:delta_alpha}
\Delta\alpha(t)=\arctan\left(\dfrac{-\dot{y}_g}{\dot{x}_g}\right),
\end{equation}
i.e. the angle between the velocity vector of the aircraft and the inertial $x_g$-axis, measured positive if the $y_g$-axis component of the velocity is negative, i.e. if the aircraft is descending. $T$ is the tension force on the tether:
\begin{equation}\label{E:Ft}
T(t)=\min\left(0,k_t\;\left(\|(x_g(t),y_g(t))\|_2-r_{M_1}\,x_1(t)\right)\right),
\end{equation}
where $k_t$ is the stiffness of the tether, assumed constant for simplicity. The saturation to 0 in Eq.~\eqref{E:Ft} accounts for the fact that the tether can only transfer force when under tension, i.e. when its length $r_{M_1}\,x_1(t)$ is smaller than the position of the aircraft relative to its attachment point on the ground. Finally, $F_l$ and $F_d$ are, respectively, the aerodynamic lift and drag forces developed by the aircraft, computed as:
\begin{equation}\label{E:FL_D}
\begin{array}{c}
F_l(t)=\frac{1}{2}\rho A C_l(\alpha(t)) \cdot \|(\dot{x}_g(t),\dot{y}_g(t))\|_2^2\\
F_d(t)=\frac{1}{2}\rho A C_{d,eq}(\alpha(t)) \cdot \|(\dot{x}_g(t),\dot{y}_g(t))\|_2^2
\end{array}
\end{equation}
where $\alpha(t)$ is the angle of attack:
\begin{equation}\label{E:alpha}
\alpha(t)=\vartheta_0+\Delta\alpha(t)+x_{9}(t).
\end{equation}

\begin{figure}
 \begin{center}
    \includegraphics[trim= 0cm 0cm 0cm 0cm,width=8cm]{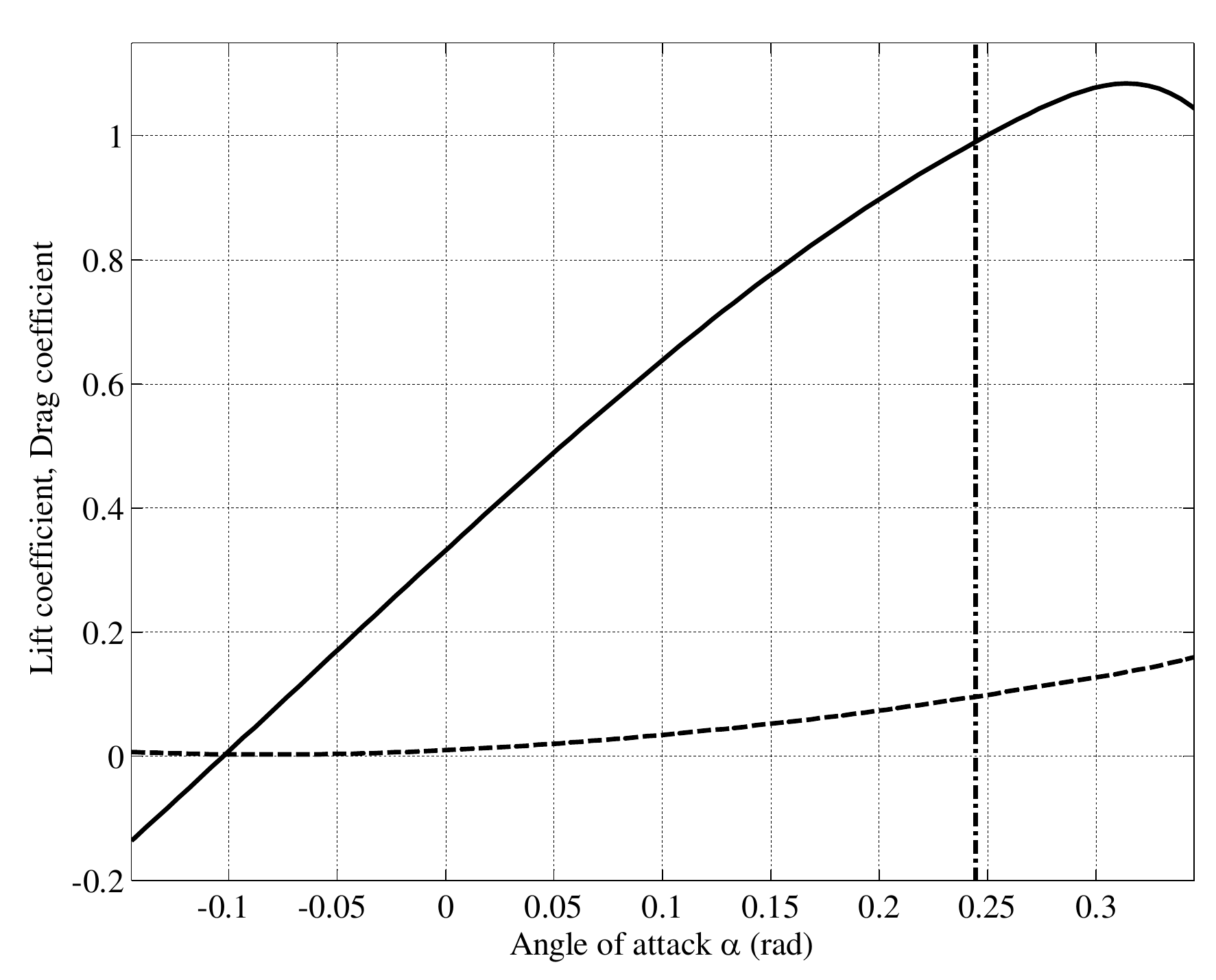}
  \caption{Lift (solid line) and drag (dashed) coefficients used in the dynamical simulation model of the take-off phase, and initial wing trimming $\vartheta_0$ (dash-dotted line).}
  \label{F:Cl_Cd_alpha}
 \end{center}
\end{figure}

The angle $\vartheta_0$ is a fixed setting for the wings' orientation, such that if the aircraft is flying horizontally (i.e. $\Delta\alpha=0$) at zero pitch angle then we have $\alpha=\vartheta_0$. The considered courses of $C_l,\,C_d$ as a function of $\alpha$ are shown in Figure \ref{F:Cl_Cd_alpha} and correspond to a finite wing with Clark-Y profile \cite{Mccor95}. The same figure also shows the chosen trimming for $\vartheta_0$.

We denote the initial state with $x^I_0$, which is required to simulate the model \eqref{E:dyn_mod_first}, i.e. $x(0)=x^I_0$. In particular, we choose the initial condition
\begin{equation}\label{E:x0_I}
x^I_0=\left[\frac{l_0}{r_{M_1}},\,0,\,0,\,0,\,x_{g,0},\,0,\,0,\,0,\,0,\,0\right]^T,
\end{equation}
meaning that the motors, the slide and the aircraft are initially at rest, a length $l_0$ of tether is reeled out and the distance of the aircraft's starting position from the attachment point of the line on the winch is equal to $x_{g,0}$, with  $x_{g,0}>l_0$ so that the tether is not exerting any force on the glider and the slide (see Eq.~\eqref{E:Ft}).

The switch between the first and the second operating modes takes place at the time instant $t^*$ defined as:
\begin{equation}\label{E:t_star}
t^*=\min\left(\tau\geq0\,:F_l(\tau)\cos(\Delta\alpha(\tau))> mg\right).
\end{equation}
Thus, $t^*$ represents the time instant when the vertical lift force developed by the glider is larger than its weight, hence obtaining a positive vertical acceleration. The initial condition $x^{II}_0$ of the model that describes the system in the second operating mode is then given by:
\begin{equation}\label{E:x0_II}
x^{II}_0=x(t^*),
\end{equation}
i.e. the state of the system in the first operating mode at the switching instant $t^*$. The model equations for the second operating mode are the following:
\begin{equation}\label{E:dyn_mod_second}
\begin{array}{rcl}
\dot{x}_1(t)&=&x_2(t)\\
\dot{x}_2(t)&=&\dfrac{1}{J_{M_1}}(r_{M_1}\,T(t)-\beta_{M_1}\,x_2(t)+u_1(t))\\
\dot{x}_3(t)&=&x_4(t)\\
\dot{x}_4(t)&=&\dfrac{1}{J_{M_2}+m_s\,r_{M_2}^2}\,(-r_{M_2}^2\,\beta_s\,x_4(t)-\beta_{M_2}\,x_4(t)\\
&
&+u_2(t))\\
\dot{x}_5(t)&=&x_6(t)\\
\dot{x}_6(t)&=&\dfrac{1}{m+m_t(t)}\,(F_l(t)\,\sin(\Delta\alpha(t))\\
&
&-F_d(t)\,\cos(\Delta\alpha(t))+\cos(x_9(t))\,u_3(t))\\
\dot{x}_7(t)&=&x_8(t)\\
\dot{x}_8(t)&=&\dfrac{1}{m+m_t(t)}\,(F_l(t)\,\cos(\Delta\alpha(t))\\
&
&+F_d(t)\,\sin(\Delta\alpha(t))\\
&
&-(m+m_t(t))\,g+\sin(x_9(t))\,u_3(t))\\
\dot{x}_9(t)&=&x_{10}(t)\\
\dot{x}_{10}(t)&=&\omega_\beta(-\Delta\alpha(t)-x_{10}(t)),
\end{array}
\end{equation}
where $m_t$ is the mass of the tether that has been reeled out:
\begin{equation}\label{E:m_t}
m_t(t)=\rho_t\,\pi\,r_t^2\,r_{M_1}\,x_1(t)
\end{equation}
with $\rho_t$ and $r_t$ being respectively the density and the radius of the tether. Regarding the last two equations in Eq.~\eqref{E:dyn_mod_second}, which describe the behavior of the pitch angle, we assume for simplicity that an active control system actuates the elevator in order to track the angle $\vartheta_{g,ref}\doteq-\Delta\alpha(t)$ with no offset, and that the resulting closed-loop dynamical behavior is given by a first-order system with time constant $\frac{1}{\omega_\beta}\,$, where $\omega_\beta$ is a constant parameter. In this way, if a steady state is attained during the ascend, the corresponding angle of attack will match the parameter $\vartheta_0$, see Eq.~\eqref{E:alpha}. Note that the pitch angle $\vartheta_g$ (i.e. $x_{9}$) affects how the thrust force $u_3$ exerted by the propeller acts on the horizontal and vertical dynamics of the aircraft, hence providing a further coupling between the pitch dynamics and the aircraft translational motion.

Eqs.~\eqref{E:dyn_mod_first}-\eqref{E:m_t} provide the hybrid model that we use to refine the results given in section \ref{S:evaluation}. However, this model cannot be simulated without first implementing suitable feedback controllers, since the open-loop behavior of the system is not stable. In the next section, we briefly describe the controllers we employ to carry out the numerical simulations.

\subsection{Control design}\label{SS:control_design}
The control objectives are different between the first and second operating mode. In the first mode, the winch motor $M_1$ has to accelerate fast enough, such that the tether tension is always zero, but avoiding at the same time that an excessive tether length is reeled-out, to limit the line sag. At the same time, the slide motor $M_2$ has to accelerate from zero to take-off speed. To achieve these goals, we employ the following proportional controllers:
\begin{equation}\label{E:control_phase1}
\begin{array}{rcl}
u_1(t)&=&K_{M_1}(\dot{x}_{g,to}-r_{M_1}\,x_2(t))\\
u_2(t)&=&K_{M_2}(\dot{x}_{g,to}-r_{M_2}\,x_2(t))\\
\end{array}
\end{equation}
where $K_{M_1},\,K_{M_2}$ are the controllers' gains, and $\dot{x}_{g,to}$ is a reference speed.

In the second operating mode, the winch motor $M_1$ shall maintain a reel-out speed that matches that of the aircraft, again to keep the tether tension at a low value. The motor $M_2$ shall brake and stop the slide. Finally, the on-board propeller shall track a desired vertical velocity $\dot{y}_{g,to}$. To obtain these goals, we employ the following proportional controllers for the motors:
\begin{equation}\label{E:control_phase2}
\begin{array}{rcl}
u_1(t)&=&K_{M_1}(\|(\dot{x}_g(t),\dot{y}_g(t))\|_2-r_{M_1}\,x_2(t))\\
u_2(t)&=&-K_{M_2}r_{M_2}\,x_2(t),\\
\end{array}
\end{equation}
while for the propeller we implement a dynamical cascade controller whose transfer function in the Laplace domain is the following
\begin{equation}\label{E:control_propeller}
C(s)\doteq\dfrac{U_3(s)}{E_{\dot{y}_g}(s)}=K_T\dfrac{\left(1+\frac{s}{w_{z,1}}\right)\left(1+\frac{s}{w_{z,2}}\right)}{s\left(1+\frac{s}{w_{p}}\right)}
\end{equation}
where $s$ is the Laplace variable, $U_3(s)$ and $E_{\dot{y}_g}(s)$ are the Laplace transforms of the propeller thrust signal $u_3(t)$ and of the tracking error $e_{\dot{y}_g}(t)\doteq\dot{y}_{g,to}-\dot{y}_g(t)$, respectively, and $K_T,\,w_{z,1},\,w_{z,2}$ and $w_{p}$ are design parameters. The need for the slightly more complex controller \eqref{E:control_propeller} for the propeller, with respect to the simple proportional gains \eqref{E:control_phase1}-\eqref{E:control_phase2} used for the motors, stems from the presence of  additional dynamics in the glider, for example due to the interaction between the pitch dynamics and the translational motion, that need to be compensated in order to avoid an oscillatory behavior of the system's response. All three inputs $u_1,\,u_2,\,u_3$ are saturated due to physical limitations of the motors:
\begin{equation}\label{E:control_sat}
\begin{array}{rcl}
-\overline{C}_{M_1}\leq&u_1(t)&\leq\overline{C}_{M_1}\\
-\overline{C}_{M_2}\leq&u_2(t)&\leq\overline{C}_{M_2}\\
0\leq&u_3(t)&\leq\overline{F_T}\\
\end{array}
\end{equation}
Finally, the described controllers are implemented in discrete time with a sampling frequency of 100$\,$Hz.

\subsection{Simulation results and discussion}\label{SS:results}
We simulate the take-off maneuver for three different aircrafts, whose effective areas matches those considered in section \ref{S:evaluation}. The model and control parameters employed for the simulations are shown in Tables \ref{T:tab_sim} and \ref{T:contr_param}, respectively. In addition, the values $\rho=1.2\,\textrm{kg/m}^{3}$,\,$g=9.81\,\textrm{m/s}^2$ and the aerodynamic coefficients shown in Figure \ref{F:Cl_Cd_alpha} have been used. The initial conditions \eqref{E:x0_I} with $l_0=2\,$m  and $x_{g,0}=0$ were used for all three aircrafts. The number and size of the propellers, required to compute the related power according to equation \eqref{heli_power}, are the same as those considered in section \ref{S:evaluation}, i.e. 2 propellers with efficiency 0.7 and 0.25$\,$m, 0.5$\,$m, 1$\,$m of diameter, respectively, for the three aircraft sizes.

\begin{table}\center
	\caption{System parameters employed to simulate the take-off maneuver. \label{T:tab_sim}}
	\begin{tabular}{|l|r|r|r|}
		\hline
		$d$ (m)& 5 & 10 & 20 \\ \hline\hline
$J_{M_1}$ (kg\,m$^2$)&1.3 & 30&490 \\\hline
$\beta_{M_1}$ (kg/s) & 0.001 &0.002 &0.003 \\\hline
$r_{M_1}$ (m)&0.2 & 0.5& 1\\\hline
$J_{M_2}$ (kg\,m$^2$)&0.03 & 0.1& 2\\\hline
$\beta_{M_2}$ (kg/s) & 0.001 & 0.002& 0.003\\\hline
$r_{M_2}$ (m)&0.1 & 0.15& 0.4\\\hline
$m_s$ (kg)& 6& 30& 120\\\hline
$m$ (kg)& 37.5& 150& 600\\\hline
$\beta_{s}$ (kg/s) & 0.1& 0.3& 1\\\hline
$k_t$ (N/m) & 1\,10$^5$& 9.1\,10$^5$& 2.5\,10$^5$\\\hline
$r_t$ (m)& 0.0025& 0.0075& 0.0125\\\hline
$\rho_t$ (kg/m$^{3}$) & 970& 970& 970\\\hline
$\omega_\beta$ (rad/s$^{3}$) & 10& 10& 10\\ \hline
$\vartheta_0$ (rad)& 0.24 & 0.24& 0.24\\ \hline
	\end{tabular}
\end{table}

\begin{table}\center
	\caption{Control parameters employed to simulate the take-off maneuver. \label{T:contr_param}}
	\begin{tabular}{|l|r|r|r|}
		\hline
		$d$ (m)& 5 & 10 & 20 \\ \hline\hline
$\dot{x}_{g,to}$ (m/s) & 30& 30& 30\\\hline
$\dot{y}_{g,to}$ (m/s) & 1& 1& 1\\\hline
$K_{M_1}$ (N\,m\,s/rad) & 3 & 20& 160\\\hline
$K_{M_2}$ (N\,m\,s/rad) & 10& 50& 200\\\hline
$K_{T}$ (N\,m\,s/rad) & 100& 150& 600\\\hline
$\omega_{p}$ (rad/s) & 16& 32& 32\\\hline
$\omega_{z,1}$ (rad/s) & 0.2& 0.2& 0.2\\\hline
$\omega_{z,2}$ (rad/s) & 1& 2& 2\\\hline
$\overline{C}_{M_1}$ (N\,m)& 750& 3000& 12000\\\hline
$\overline{C}_{M_2}$ (N\,m)& 48& 290& 3500\\\hline
$\overline{F_T}$ (N)& 80& 350& 600\\\hline
	\end{tabular}
\end{table}

\begin{figure}[!htb]
 \begin{center}
  \includegraphics[trim= 0cm 0cm 0cm 0cm,width=8cm]{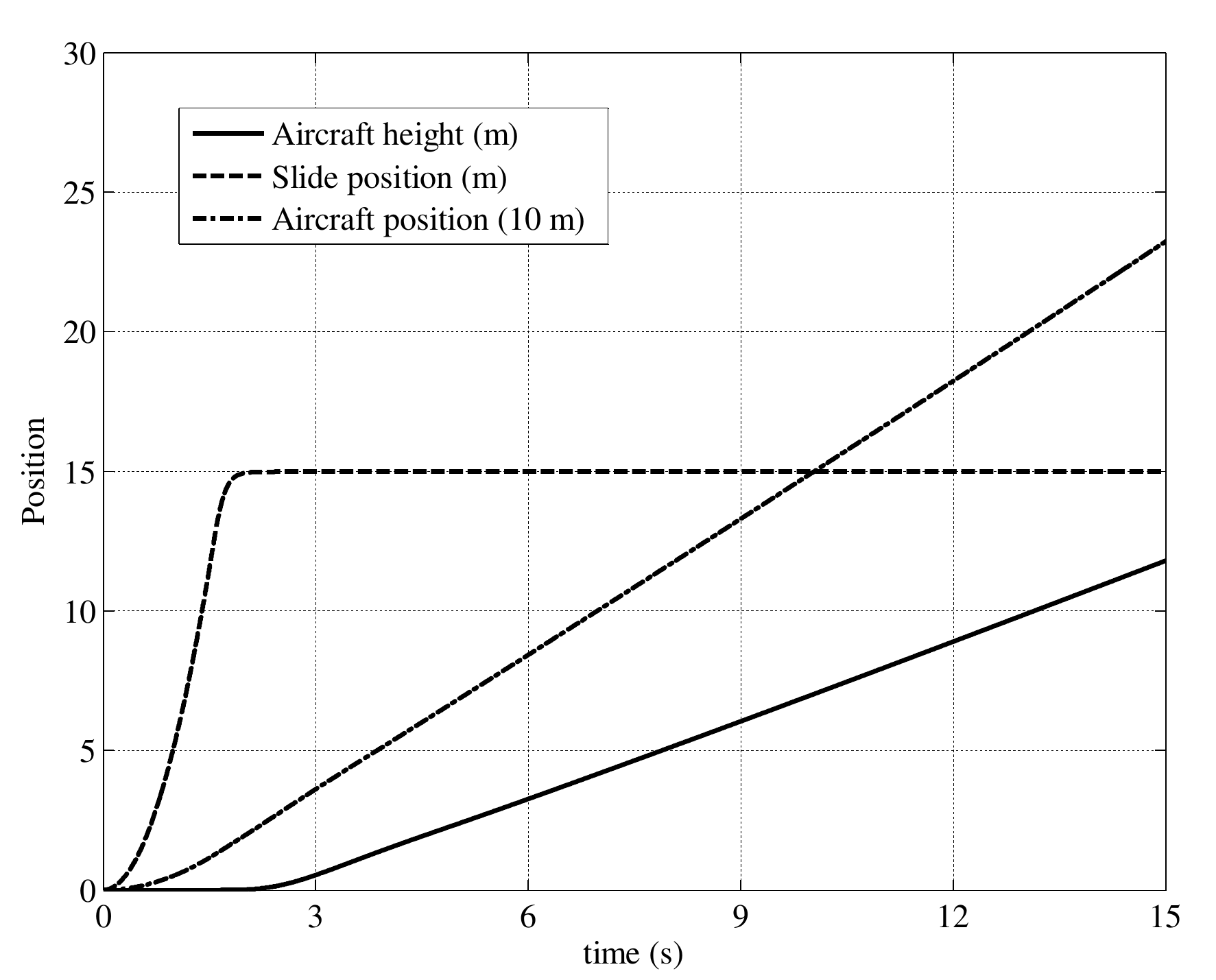}
  \caption{Simulation results with the $10$-m-wingspan aircraft. Courses of the aircraft height, slide position and aircraft distance from the ground station (divided by 10 for the sake of clarity).}
  \label{F:sim_pos}
 \end{center}
\end{figure}

\begin{figure}[!htb]
 \begin{center}
  \includegraphics[trim= 0cm 0cm 0cm 0cm,width=8cm]{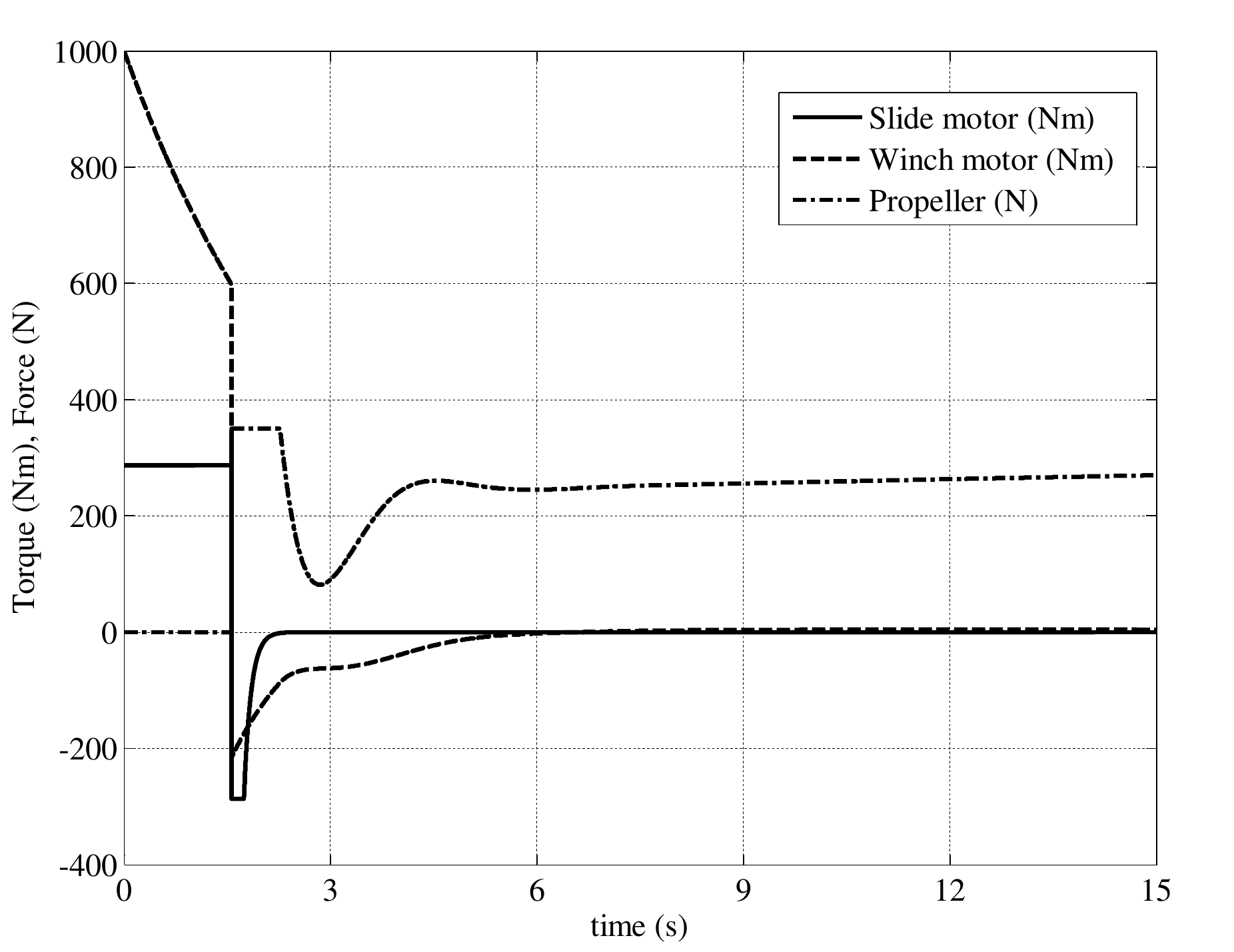}
  \caption{Simulation results with the $10$-m-wingspan aircraft. Courses of the motor torques and of the propeller thrust.}
  \label{F:torque_force}
 \end{center}
\end{figure}

\begin{figure}[!htb]
 \begin{center}
  \includegraphics[trim= 0cm 0cm 0cm 0cm,width=8cm]{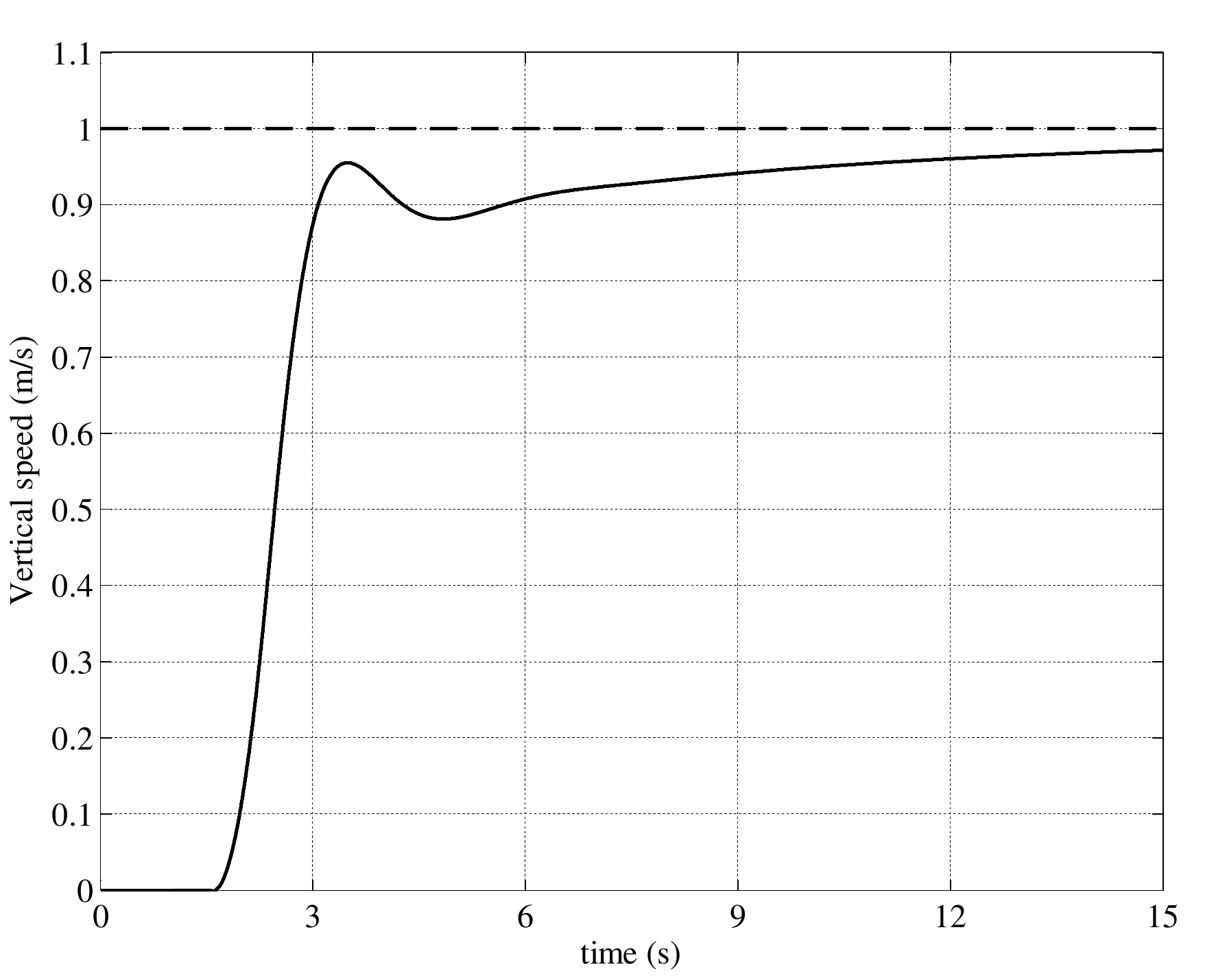}
  \caption{Simulation results with the $10$-m-wingspan aircraft. Course of the vertical speed of the aircraft (solid) and the target value (dashed).}
  \label{F:vertical_speed}
 \end{center}
\end{figure}

\begin{figure}[!htb]
 \begin{center}
  \includegraphics[trim= 0cm 0cm 0cm 0cm,width=8cm]{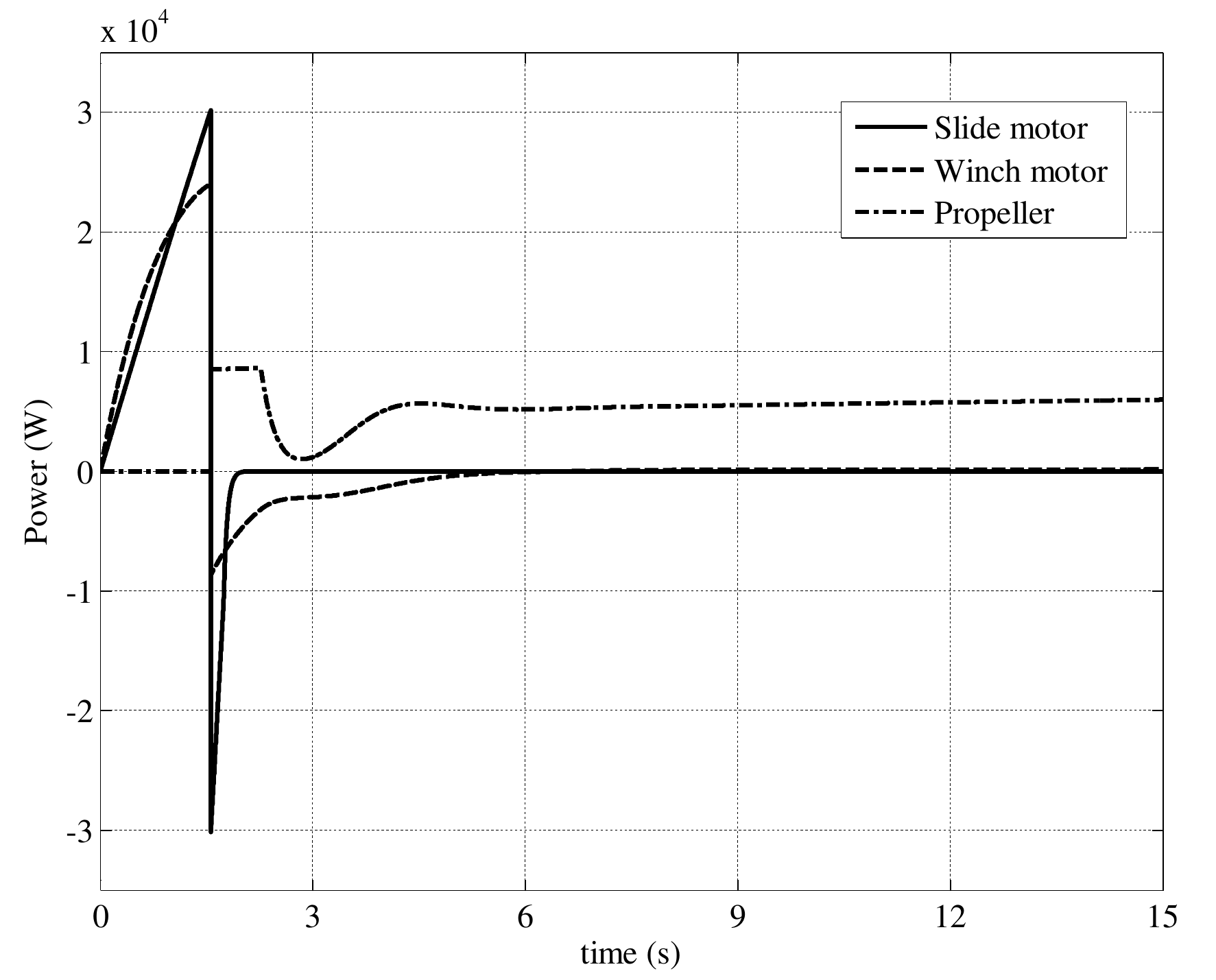}
  \caption{Simulation results with the 10-m wingspan aircraft. Courses of the motors' and propeller's power.}
  \label{F:power}
 \end{center}
\end{figure}

Examples of simulation results for the aircraft with $d=10\,$m are shown in Figures \ref{F:sim_pos}-\ref{F:power}. In Figure \ref{F:sim_pos}, it can be noted that the total travel distance of the slide is equal to 15 m, and that the aircraft starts the ascend after 12.4 m, i.e. when the take-off speed of 15.7 m/s has been reached. As shown in Figure \ref{F:torque_force}, the motor $M_2$ exploits the full rated torque to accelerate and then to brake the slide, while $M_1$ employs a relatively small fraction of its available torque for the acceleration and then settles to a constant torque corresponding to the viscous friction at the aircraft's velocity. We remark that the power required to accelerate the drum, although substantial, does not give rise to additional costs, since the machine $M_1$ is already present and the power required for take-off is a small fraction of the one that occurs during power generation. The propeller is engaged only after take-off and, after a short transient, it settles to a steady value sufficient to achieve the desired vertical velocity. The behavior of the latter quantity as compared with its reference is reported in Figure \ref{F:vertical_speed}. As shown in Figure \ref{F:power}, the peak power for the motors is reached at the instant when the aircraft takes off. The results obtained with the other two aircrafts ($d=5$ and 20 m) are qualitatively similar to those shown in Figures \ref{F:sim_pos}-\ref{F:power}. In all cases, the total travel distance of the slide was about 15 m.

\begin{table*}\center
	\caption{Comparison between the power values provided by the simplified equations and those provided by the numerical simulations. The percentages in brackets refer to the peak mechanical power of the generator with 15 m/s wind speed.\label{T:sim_results}}
	\begin{tabular}{|l|c|c|c|}
		\hline
		Wingspan (m)& 5 & 10 & 20 \\ \hline\hline
Ground motor (kW) - simple equation& 8 (11\%)& 31 (10\%)& 124 (10\%)\\\hline
Ground motor (kW) - simulation& 8 (11\%)& 30 (10\%)& 140(11\%)\\\hline\hline
Propeller (kW) - simple equation& 2 (3\%)& 9 (3\%)& 37 (3\%)\\\hline
Propeller (kW) - simulation& 3 (4\%)&13 (4\%)& 50 (4\%)\\\hline
	\end{tabular}
\end{table*}

Table \ref{T:sim_results} shows a comparison between the power figures obtained from the simplified analysis of section \ref{S:evaluation} and those obtained with the simulations. The  values of power required on the ground are very well matching, hence confirming the outcome of our simplified analysis. The larger simulated values for the required on-board power, with respect to the simplified analysis, are due to the inertia of the aircraft, which plays a role in the transient from zero vertical speed to the target one (see Figure \ref{F:torque_force}), and due to its pitch, which has the effect of decreasing the thrust in horizontal direction and adding a braking contribution from the lift force projected onto the $x_g-$axis. Again, notwithstanding these effects, the on-board power required for the ascend appears to be a reasonable fraction of the system's power. Moreover, we did not carry out any optimization neither of the design parameters nor of the controllers, which can still be adjusted in order to achieve different tradeoffs between peak power consumption and velocity of the transient from zero to the target vertical speed.

\section{Conclusions}\label{S:conclusions}
We presented an analysis of different concepts for the take-off phase of AWE systems based on rigid wings and ground-level power conversion, by means of basic equations. Based on the analysis, we concluded that a linear take-off maneuver with a ground acceleration phase and on-board propellers is the most viable approach. We refined the analysis of this maneuver by means of numerical simulations with a hybrid dynamical model. The simulation results predict slightly larger on-board power values than the simplified analysis, but still they are small compared to the total power of the generator. This indicates that the take-off equipment constitutes a rather small cost fraction of the total system costs. At the same time, the required land occupation appears to be reasonable. These outcomes confirm the technical and economic feasibility of this take-off technique. Further studies will be devoted to a deeper analysis of the approach and to the study of the landing maneuver, both with finer dynamical models, also accounting also for wind turbulence, and with experimental activities.





\bibliographystyle{plain}

\end{document}